\documentclass[11pt,a4paper]{article}

\usepackage{amsmath,amsfonts,amssymb,latexsym,graphics,epsfig,url}
\usepackage{xcolor}
\usepackage{amsthm}
\usepackage[english]{babel}
\usepackage{mathdots}
\usepackage{graphicx}
\usepackage{soul} 
\usepackage[utf8]{inputenc}
\usepackage{diagbox}
\usepackage[makeroom]{cancel}
\usepackage{multicol}
\usepackage{float}              
\newtheorem{theorem}{Theorem}[section]
\newtheorem{proposition}[theorem]{Proposition}

\newtheorem{conjecture}[theorem]{Conjecture}
\newtheorem{example}[theorem]{Example}

\newtheorem{definition}[theorem]{Definition}

\DeclareMathOperator{\spec}{sp}

\def\Z{\ns Z}

\def\vecu{\mbox{\boldmath $u$}}
\def\vecv{\mbox{\boldmath $v$}}

\def\vec0{\mbox{\boldmath $0$}}

\def\A{\mbox{\boldmath $A$}}
\def\B{\mbox{\boldmath $B$}}

\def\G{\Gamma}

\def\M{\mbox{\boldmath $M$}}

\def\Z{\ns{Z}}

\def\1{\mbox{\boldmath $1$}}

\def\G{\Gamma}

\def\Re{\mathbb R}
\def\Z{\mathbb Z}

\begin{document}
	
\title{Structural and Spectral Properties of Chordal Ring, Multi-ring and Mixed Graphs
\thanks{This research has been supported by
AGAUR from the Catalan Government under project 2021SGR00434 and MICINN from the Spanish Government under project PID2020-115442RB-I00.
M. A. Fiol's research was also supported by a grant from the  Universitat Polit\`ecnica de Catalunya with references AGRUPS-2022 and AGRUPS-2023.}
}
	\author{M. A. Reyes$^a$, C. Dalf\'o$^a$, M. A. Fiol$^b$\\
		\\
		{\small $^a$Dept. de Matem\`atica, Universitat de Lleida, 08700 Igualada (Barcelona), Catalonia}\\
		{\small {\tt \{monicaandrea.reyes,cristina.dalfo\}@udl.cat}}\\
		{\small $^{b}$Dept. de Matem\`atiques, Universitat Polit\`ecnica de Catalunya, 08034  Barcelona, Catalonia} \\
		{\small Barcelona Graduate School of Mathematics} \\
		{\small  Institut de Matem\`atiques de la UPC-BarcelonaTech (IMTech)}\\
		{\small {\tt miguel.angel.fiol@upc.edu} }\\
	}

\date{}
\maketitle
	
\begin{abstract}
The chordal ring (CR) graphs are a well-known family of graphs used to model some interconnection networks for computer systems in which all nodes are in a cycle. 
Generalizing the CR graphs, in this paper, we introduce the families of chordal multi-ring (CMR), chordal ring mixed (CRM), and chordal multi-ring mixed (CMRM) graphs. In the case of mixed graphs, we can have edges (without direction) and arcs (with direction). The chordal ring and chordal ring mixed graphs are bipartite and 3-regular. They consist of a number $r$ (for $r\geq 1$) of (undirected or directed) cycles with some edges (the chords) joining them. In particular, for CMR, when $r=1$, that is, with only one undirected cycle, we obtain the known families of chordal ring graphs.
Here, we use plane tessellations to represent our chordal multi-ring graphs. This allows us to obtain their maximum number of vertices for every given diameter. Besides, we computationally obtain their minimum diameter for any value of the number of vertices.
Moreover, when seen as a lift graph (also called voltage graph) of a base graph on Abelian groups, we obtain closed formulas for the spectrum, that is, the eigenvalue multi-set of its adjacency matrix.
\end{abstract}
	
	

\section{Introduction}

Interconnection networks play a crucial role in the design of distributed computer systems, often represented by graphs where vertices symbolize nodes or processing elements, and edges denote communication links between them. See, for instance, Deo \cite{d74}, and Bermond et al. \cite{bdq86}. The topology of these networks is a key factor influencing communication delay, throughput, and message routing efficiency. Understanding graph theory concepts such as adjacency, degree, diameter, and mean distance is essential in analyzing and optimizing these networks. Exploring topologies like rings, meshes, and hypercubes can offer efficient and reliable interconnection solutions for multi-computer systems. 
One of the simplest topologies for interconnection networks is the undirected or directed graphs in
which each node is connected to two others, making up a bidirectional loop or cycle. The main drawbacks of cycles are their poor
reliability (any link or processor failure disconnects the network) and low
performance (some messages must travel along half or the whole ring to reach their destination).
To overcome these problems, the cyclic topology is improved by the chordal ring graphs, which consist of a cycle with some additional links between nodes. 
 Arden and Lee [3] proposed the chordal ring graphs for efficient and reliable multi-(micro)computer interconnection networks, which are real-world examples. 
 If only one link is added to each node, the corresponding graphs are $3$-regular. This is the case of the so-called `chordal ring
networks', see an example in Figure \ref{fig:qcr} (left).
Chordal rings were first introduced by Coxeter \cite{c50}. Since then, the structure of these graphs has been extensively studied. For example, Arden and Lee \cite{al81} studied the problem of the maximization of the number of nodes for a given diameter, and Yebra et al. \cite{yfma85} found a relationship between certain types of plane tessellations, where the vertices are associated with regular polygons and chordal ring graphs (see also Morillo et al. \cite{mcf87}, 
Dalf\'o et al. \cite{deeflmt24}). This geometrical representation
characterizes the graph and facilitates the study of some of its parameters, particularly those with distance-related parameters.
In this paper, we show that chordal ring graphs are one of four closely related families of graphs constituted by (undirected or directed) cycles with chords. These are the chordal ring and multi-ring graphs ($CR$ and $CMR$) and the chordal ring and multi-ring mixed graphs ($CRM$ and $CMRM$). 

With the chordal graphs, one finds a smaller diameter for any value of the number of vertices and a greater number of vertices for every given diameter, both with respect to the cycle graphs.
Our research gap is to improve these results with the multi-ring chordal graphs and the ring and multi-ring mixed graphs. 
This unified approach allows us to represent all of them as quotient graphs of a `running bond' infinite graph-pattern. This name refers to a way of placing blocks used by bricklayers; see  Reid \cite{r04} and  Figures \ref{fig:diametre5-undirected} and \ref{fig:diameter5-mixed}. As shown on the left in the same figures, such quotient graphs are obtained when the infinite graph is taken modulo some integral matrix $\M$,  whose rows are the translation vectors defining the periodicity. More precisely, the vertices of the infinite graph are identified with integral vectors, and two vectors $\vecu$ and $\vecv$ represent the same vertex of the quotient graph if and only if $\vecu\equiv \vecv \mod \M$. That is, $\vecu-\vecv$ belongs to the lattice generated by the rows of $\M$, and $\vecu-\vecv\in \Z^2\M$, see Fiol \cite{f87}.
This approach is valid for graphs with edges, that is, bidirectional links. It is also valid for mixed graphs with edges (bidirectional links) and arcs (one-directional links). Moreover, it can be applied to any value of the number of cycles (with at least one cycle) and any number of chords (with at least one chord). The chords are the additional edges that allow us to obtain a chordal-
ring graph from a cycle.

\begin{figure}[ht]
    \centering
    \includegraphics[width=10cm]{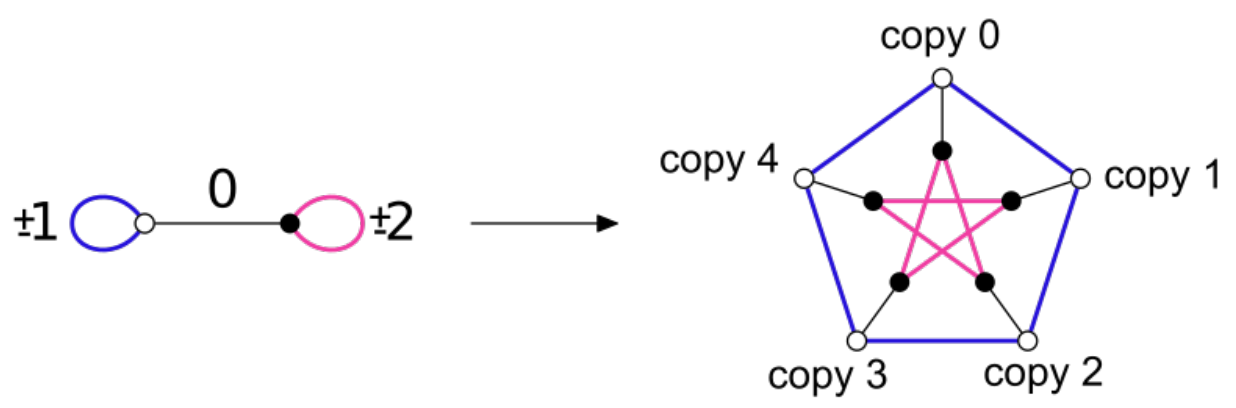}
    \caption{The base graph with a voltage assignment on the arcs (left), when acting on the group $\mathbb{Z}_5$, gives the Petersen graph (right).}
    \label{fig:petersen}
\end{figure}

We show how to construct the chordal ring families from graphs admitting a voltage assignment, that is, with `weights' on the arcs. With this respect,
Gross introduced the following concepts in \cite{g74}. Let $G$ be a group. An ({\em ordinary\/}) {\em voltage assignment} on the (di)graph $\Gamma=(V,E)$ is a mapping $\alpha: E\to G$ with the property that $\alpha(a^{-1})=(\alpha(a))^{-1}$ for every arc $a\in E$. Thus, a voltage assigns an element $g\in G$ to each arc of the graph so that a pair of mutually reverse arcs $a$ and $a^{-1}$, forming an undirected edge, receive mutually inverse elements $g$ and $g^{-1}$. The graph $\Gamma$ and the voltage assignment $\alpha$ determine a new graph $\Gamma^{\alpha}$, called the {\em lift} of $\Gamma$, which is defined as follows. The vertex set of the lift is the Cartesian product $V^{\alpha}=V\times G$. Moreover, for every arc $a\in E$ from a vertex $u$ to a vertex $v$ for $u,v\in V$ (possibly, $u=v$) in $\Gamma$, and for every $g\in G$, there is an arc $(a,g)\in E^{\alpha}$ from the vertex $(u,g)\in V^{\alpha}$ to the vertex $(v,g\alpha(a))\in V^{\alpha}$. Let us show the example illustrated in Figure \ref{fig:petersen}. In this example, when the base graph with a voltage assignment on the arcs (drawn on the left) is applied to the group $\mathbb{Z}_5$, we obtain the Petersen graph (drawn on the right). Since the group is $\mathbb{Z}_5$, we have 5 copies (numbered from 0 to 4) of the vertices of the base graph. So, as the voltage of the pink edge is $2$, we have to join the black vertices by adding $2 \mod 5$ to each copy. We do the same for the other edges, each with its corresponding voltage. 

One of this paper's main contributions is finding the spectra of the different families of chordal ring graphs and mixed graphs. We recall the conjecture by Haemers \cite{h16} that states that almost all graphs are determined by their spectra. More precisely, among all non-isomorphic graphs on, at most, $n$ vertices, the fraction of graphs that are not determined by their spectra goes to 1 when $n$ goes to infinity.
With this aim, we use a quotient-like matrix, introduced by Dalf\'o et al. \cite{dfmr17}, that fully represents a lifted digraph. The main advantage of this approach is that such a matrix has a size equal to the order of the base digraph.
For a digraph $\G=(V,E)$ with voltage assignment $\alpha$, if we deal with the case when the group $G$ of the voltage assignments is cyclic (that is,  $G=\Z_k=\{0,1,\ldots,k-1\}$), then its
\emph{polynomial matrix} $\B(z)$ is a square matrix indexed by the vertices of $\G$. The matrix elements are complex polynomials in
the quotient ring $\Re_{k-1}[z]=\Re[z]/(z^k)$, where $(z^k)$ is the ideal generated by the polynomial
$z^k$. More precisely, each entry of $\B(z)$ is fully represented by a polynomial of degree at most $k-1$, say $(\B(z))_{uv}=p_{uv}(z)=\alpha_0+\alpha_1 z+\cdots +\alpha_{k-1}z^{k-1}$, where
$$
\alpha_i=\left\{
\begin{array}{ll}
1 & \mbox{if  $uv\in E$ and $\alpha(uv)=i$,}\\
0 & \mbox{otherwise,}
\end{array}
\right.
$$
for $i=0,\ldots,k-1$.

This paper is structured as follows. In the next section, we present our approach to chordal ring graphs and obtain their spectra when seeing them as lift graphs. 
In Section \ref{sec:CMR}, we introduce the chordal multi-ring graphs to obtain graphs with an optimal number of vertices for a given diameter and find their spectra through voltage graphs. Similar results are obtained in Section \ref{sec:CRM} for the chordal ring mixed graphs (in which the cycle is directed and the chords are undirected) and, finally, in Section \ref{sec:CMRM} for the chordal multi-ring mixed graphs. Our approach allows us to
determine the spectra in all the cases.

\section{Chordal ring graphs}
\label{sec:CR}

\begin{figure}[ht]
    \centering
    \includegraphics[width=14cm]{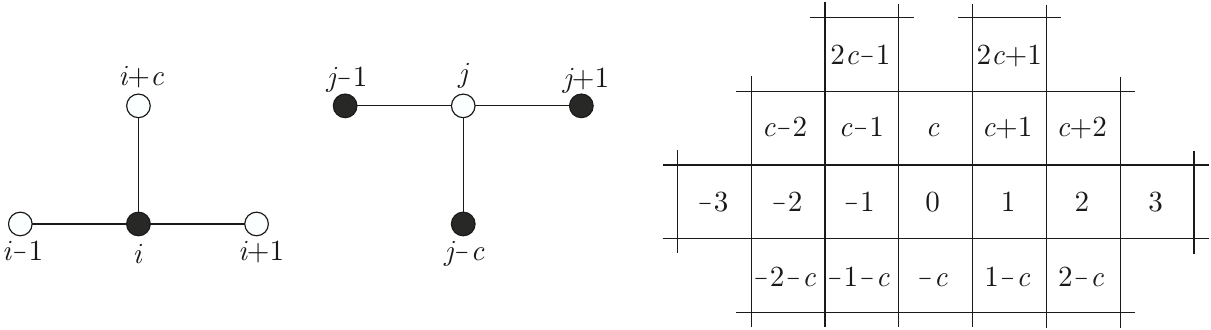}
    \caption{Adjacencies (for $i$ even and $j$ odd) and planar pattern of the vertices in the chordal ring graph $CR(N,c)$ (with vertices at a distance at most three from vertex $0$).}
    \label{fig:planar-pattern}
\end{figure}

The \textit{chordal ring graph} $CR(N,c)$ has an even number of vertices $(N=2n)$ labeled with the integers $\{0,1,\ldots,N-1\}$, and each {\em even} vertex $i$ is connected to the vertices $i\pm 1\ \mod N$ and $i+c\ \mod N$ for some odd integer $c$. Consequently, each {\em odd} vertex $j$ is connected to the vertices $j\pm 1\ \mod N$ and $j-c\ \mod N$. See Figure \ref{fig:planar-pattern} on the left. Therefore, we have a ring structure with additional links called {\em chords}. An example is the Heawood graph with a diameter 3, isomorphic to $CR(14,5)$, and it is known to be a $(3,6)$-cage, see Figure \ref{fig:qcr} (left).
The chordal ring graphs are bipartite and vertex-symmetric; that is, they have an automorphism group that acts transitively on the vertices. Recall that a group of automorphisms is an algebraic structure that defines the symmetries in the graph. More precisely, an automorphism of a graph $G = (V, E)$ is a permutation $\sigma$ of the vertex set $V$, such that the pair of vertices $(u, v)$ forms an edge if and only if the pair $(\sigma(u), \sigma(v))$ also forms an edge.  
In fact, this property is shared by the other three families studied in this paper: the multi-ring graphs $(CMR)$, and the chordal ring and multi-ring mixed graphs ($CRM$ and $CMRM$).
More details about the symmetries of the chordal ring graphs are in the following result.

\begin{proposition}
\begin{itemize}
\item[$(i)$] 
The chordal ring graph $CR(N,c)$ is isomorphic to $CR(N,-c)$.
\item[$(ii)$]
The automorphism group of $CR(N,c)$, with $N=2n$, contains the dihedral group $D_n$ with $N=2n$ elements and  presentation
\begin{equation}
D_n=\langle \sigma,\tau\, |\, \sigma^n=\tau^2=(\sigma\tau)^2=e\rangle
\label{Dn}
\end{equation}
where $e$ denotes the identity element.
\end{itemize}
\end{proposition}
\begin{proof}
$(i)$ Note that, in $CR(N,-c)$, each vertex $i$ is adjacent to $i-c$ if $i$ is even and $i+c$ if $i$ is odd. Then, let us prove that the mapping $\gamma:i\mapsto i+1$ is an isomorphism from  
$CR(N,c)$ to $CR(N,-c)$.
\begin{itemize}
\item
In $CR(N,c)$, vertex $i$ is adjacent to $i\pm 1$ whereas, in $CR(N,-c)$, $\gamma(i)=i-1$ is adjacent to $\gamma(i\pm 1)=\{i,i-2\}$.
\item 
 In $CR(N,c)$, vertex $i$ (even) is adjacent to $i+c$ whereas, in $CR(N,-c)$, $\gamma(i)=i-1$ (odd) is adjacent to $\gamma(i+c)=i-1+c$.
 \item 
 In $CR(N,c)$, vertex $i$ (odd) is adjacent to $i-c$ whereas, in $CR(N,-c)$, $\gamma(i)=i-1$ (even) is adjacent to $\gamma(i-c)=i-1-c$.
 \end{itemize}
$(ii)$ 
 Let us consider the mappings $\sigma: i\mapsto i+2$ and $\tau:i\mapsto -i-1$ from $CR(N,c)$ to itself.
Then, let us first check the defining relations in \eqref{Dn}: We have $\sigma^n(i)=i+2n=i$, $\tau^2(i)=\tau(-i-1)=-(-1-i)-1=i$, and 
$$
(\sigma\tau)^2(i)=\sigma\tau\sigma(-i-1)=\sigma\tau(-i+1)=\sigma(i-2)=i.
$$
In order to prove that $\sigma$ and $\tau$ are automorphisms,  let us represent the vertices adjacent to $i$ as $G_{\pm}(i)=\{i+1,i-1\}$ and $G_c(i)=i+c$ for $i$ even, and $G_c(i)=i-c$ for $i$ odd. Then
\begin{itemize}
\item For every vertex $i$,
$$
\sigma(G_{\pm}(i))=\sigma(i\pm 1)=\{i+3,i+1\}=G_{\pm}(i+2)=G(\sigma(i)).
$$
\item
If $i$ is even,
$$
\sigma(G_{c}(i))=\sigma(i+c)=i+2+c=G_{c}(i+2)=G_c(\sigma(i)).
$$
\item
If $i$ is odd,
$$
\sigma(G_{c}(i))=\sigma(i-c)=i+2-c=G_{c}(i+2)=G_c(\sigma(i)).
$$
\end{itemize}
Thus, in all the cases, $\sigma$ and $\tau$ commute with $G_{\pm}$ and $G_c$ proving that they are automorphisms of $CR(N,c)$.
\end{proof}

From the above result, if $i$ and $j$ are vertices with the same parity, then 
$\sigma^{\frac{j-i}{2}}(i)=j$. Otherwise, if $i$ and $j$ have different parity,
$\sigma^{\frac{i+j+1}{2}}\tau(i)=\sigma^{\frac{i+j+1}{2}}(-i-1)=j$. This shows that, as it was commented, there is always an automorphism mapping a vertex $i$ to a vertex $j$.

As commented in the Introduction, chordal ring graphs can be represented as congruent tiles that tessellate the plane periodically.
 More precisely, if each vertex of $CR(N,c)$ is represented by a numbered $\mod N$ unit square, the vertices reached at a distance $0,1,2,\ldots$ from any given vertex can be arranged in a planar pattern,  as shown in Figure \ref{fig:planar-pattern} starting from vertex $0$.

Since there are $3\ell$ vertices at a distance $\ell(>0)$ from vertex $0$, and the graph is bipartite, it was shown (in Yebra et al. \cite{yfma85}, and Morillo et al. \cite{mcf87}) that the maximum number $N(k)$ of vertices of a chordal ring with a diameter $k$ is 
\begin{equation}
N(k)=\left\{ 
\begin{array}{ll}
\frac{3k^2+1}{2} & \mbox{for $k$ odd,}\\[.1cm]
\frac{3k^2}{2} & \mbox{for $k$ even.}
\end{array}
\right.
\label{boundCR}
\end{equation}
Moreover, it was shown that such a maximum can be attained when $k$ is odd and cannot be attained when $k>2$ is even. Consequently, the following conjecture was raised in the same papers (see also Comellas and Hell \cite{ch03}).

\begin{conjecture}
\label{conjec-1}
The maximum number $N$ of vertices of a chordal ring graph with an even diameter $k>2$ is $N=\frac{3k^2}{2}-k$.
\end{conjecture}

\begin{figure}[t]
    \centering
    \includegraphics[width=15cm]{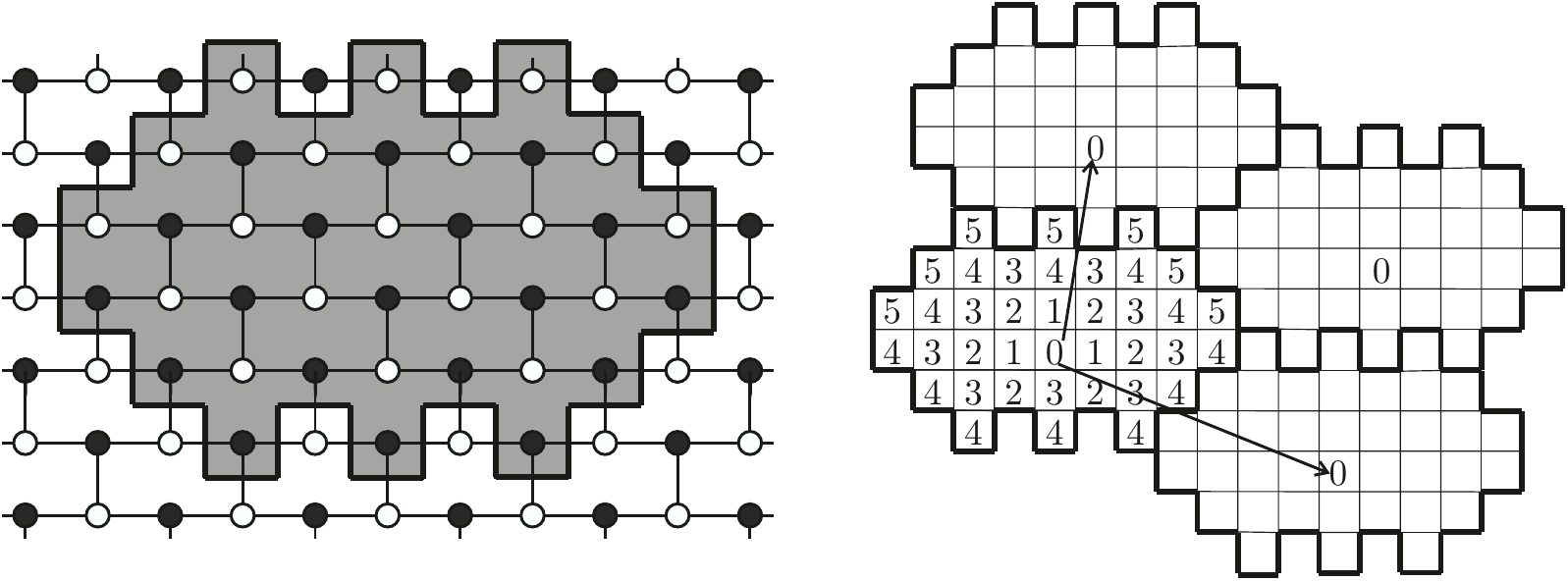}
    \caption{The optimal tiles for the chordal ring graphs with a diameter $k=5$, with translation vectors $(1,5)^{\top}$ and 
    $(7,-3)^{\top}$.}
    \label{fig:diametre5-undirected}
\end{figure}

The method used in \cite{yfma85, mcf87} to obtain chordal rings with maximum number of vertices $(3k^2+1)/2$ (for $k$ odd) and $3k^2/2-k$ (for $k$ even) consists of the following two steps: 
\begin{itemize}
\item
[{\bf $(1)$}] 
Using the planar pattern, find the `optimal tiles' (that is, containing the maximum number of vertices for a given diameter $k$), and check that they periodically tessellate the plane (see Figure \ref{fig:diametre5-undirected} for $k=5$); 
\item 
[{\bf $(2)$}]
From the two basic translation
vectors of the periodic tiling (generating the lattice of the positions of the vertices with label $0$), solve a linear system of equations to find the chord $c$.
\end{itemize}

In Table \ref{tab:tabla-cr}, we show the minimum diameter $k$ and chord $c$ for each number of vertices $N\le 528$ of a chordal ring graph $CR(N,k)$. The cases in which we get the maximum number of vertices for a given diameter are in boldface. Notice that, for an odd diameter $k$, such a maximum is as expected, whereas, for an even diameter $k$, the maximum supports Conjecture \ref{conjec-1}.

Other properties of the chordal ring graphs were studied by
Barri\`ere et al. \cite{bcm01} (gossiping),
Barri\`ere et al. \cite{bfsz00} (fault-tolerant routing), and Zimmerman and Esfahanian \cite{ze92} (fault-tolerance).

\subsection{Chordal ring graphs as lifts}
\label{sec:CR-lifts}

\begin{figure}[ht]
    \centering
   \includegraphics[width=10cm]{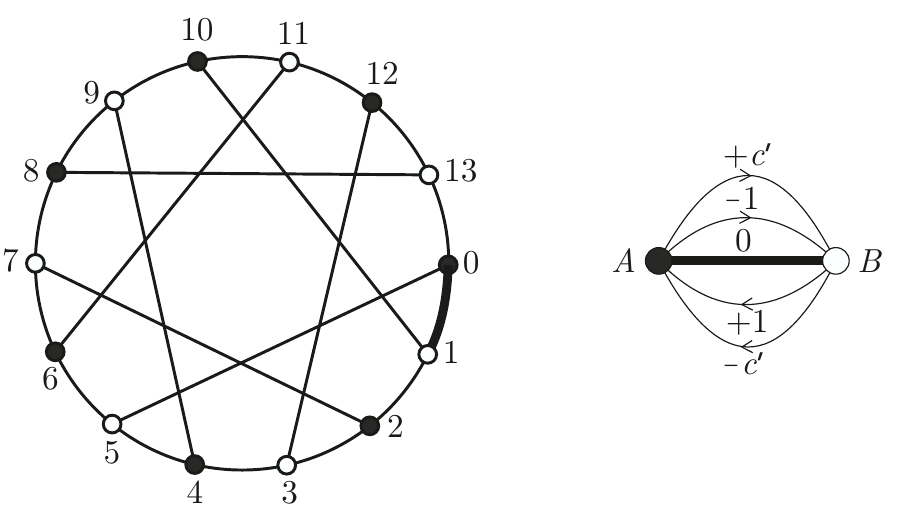}
    \caption{Left: The Heawood graph with a diameter 3, which is isomorphic to the chordal ring $CR(14,5)$. Right: The base graph of a chordal ring graph $CR(N,c)$ on the group $\Z_{N/2}$, where $c'=(c-1)/2$. (The black and white vertices stand for the even and odd vertices, respectively). The thick lines represent the edges with voltage 0 in the base graph and in the first copy of the lift graph.}
    \label{fig:qcr}
\end{figure}

The chordal ring graph $CR(N,c)$, with $N=2n$, can be seen as a lift of a base graph on the group $\Z_n$, which is represented in Figure \ref{fig:qcr} (right). This allows us to derive a closed formula giving all the eigenvalues of $CR(N,c)$.

\begin{proposition}
Given integers $N=2n$ and $c$ (odd), the eigenvalues of the chordal ring graph $CR(N,c)$ are
\begin{equation}
\lambda(r)_{1,2}=\pm\sqrt{4\cos^2\left(\frac{r\pi}{n}\right)+4\cos^2\left(\frac{(c-1)r\pi}{2n}\right)+4\cos^2\left(\frac{(c+1)r\pi}{2n}\right)-3} 
\label{specCR(N,c)}
\end{equation}
for $r=0,1,\ldots,n-1$.
\end{proposition}

\begin{proof}
 The polynomial matrix of the base graph in Figure \ref{fig:qcr} is
$$
		\B(z)=
		\left(
		\begin{array}{cc}
		0 & 1+\frac{1}{z}+z^{c'}  \\
		1+z+\frac{1}{z^{c'}} & 0 
		\end{array}
		\right),
$$
where $c'=\frac{c-1}{2}$. Then, the eigenvalues of $CR(N,c)$ can be obtained as the eigenvalues of $\B(z)$ for every $z=\zeta^r$ with $\zeta=e^{\frac{i2\pi}{n}}$, for $n=N/2$ and $r=0,1,\ldots,n-1$. With  $\tau=1+\frac{1}{z}+z^{c'}$, such eigenvalues are 
\begin{align*}
\lambda(r)_{1,2}&=\pm\sqrt{\tau\overline{\tau}}=\pm|\tau|\\
 &=\pm\sqrt{\left[1+\cos\left(\frac{r2\pi}{n}\right)+\cos\left(\frac{c'r2\pi}{n}\right)\right]^2
 +\left[-\sin\left(\frac{r2\pi}{n}\right)+\sin\left(\frac{c'r2\pi}{n}\right)\right]^2}
\end{align*}
and, operating with $c'=\frac{c-1}{2}$, we obtain \eqref{specCR(N,c)}.
\end{proof}

\begin{example}
		In the case of the graph  $CR(20,5)$, we list its eigenvalues \eqref{specCR(N,c)} for every $r=0,1,\ldots,9$ in Table \ref{tab:ev-cr}.
  \begin{table}[ht]
  \caption{All the eigenvalues of the matrices $\B(\zeta^{r})$, which yield the eigenvalues of the chordal ring graph $CR(20,5)$.}
		\label{tab:ev-cr}
		\begin{center}
	\begin{tabular}{|c|c|}
	\hline
	$r$ & $\lambda(r)_{1,2}$ \\
	\hline\hline
	$0$ &  $\pm 3$ \\
     \hline
     \rule{0pt}{3.5ex} $1$ &  $\pm\sqrt{3+2\cos\left(\frac{\pi}{5}\right)}\approx \pm 2.149$\\[1.5ex]
     \hline
     \rule{0pt}{3ex} $2$ &  $\pm 2\cos\left(\frac{\pi}{5}\right)-1\approx \pm 0.6164$ \\[1ex]
      \hline
     \rule{0pt}{3.5ex} $3$ &  $\pm\sqrt{3-2\cos\left(\frac{\pi}{5}\right)}\approx \pm 1.543$ \\[1.5ex]
      \hline
     \rule{0pt}{3ex} $4$ &  $\pm\frac{1}{2}(1+\sqrt{5})\approx \pm 1.618$ \\[1ex]
      \hline
     $5$ &  $\pm 1$ \\
      \hline
     \rule{0pt}{3ex} $6$ &  $\pm\frac{1}{2}(1+\sqrt{5})\approx \pm 1.618$ \\[1ex]
      \hline
      \rule{0pt}{3.5ex} $7$ &  $\pm\sqrt{3+2\cos\left(\frac{\pi}{5}\right)}\approx \pm 1.543$ \\[1.5ex]
       \hline
     \rule{0pt}{3ex} $8$ &  $\pm 2\cos\left(\frac{\pi}{5}\right)-1\approx \pm 0,6164$ \\[1ex]
     \hline
     \rule{0pt}{3.5ex} $9$ &   $\pm\sqrt{3+2\cos\left(\frac{\pi}{5}\right)}\approx \pm 2.149$\\[1.5ex]
					\hline
				\end{tabular}
		\end{center}
		\end{table}
	\end{example}


\section{Chordal multi-ring graphs}
\label{sec:CMR}

To generalize chordal ring graphs, we introduce the family of chordal multi-ring graphs.

\begin{definition}
\label{def:CMR}
Given positive integers $m$, $n$
(even), and $c(>1)$ (odd), the chordal 
$m$-ring graph $CMR(m,n,c)$ has vertices labeled with the elements of the Abelian group $\Z_m\times \Z_n$, and edges $(\alpha,i)\sim (\alpha,i\pm 1)$ for every $\alpha\in \Z_m$ and $i\in \Z_n$, and $(\alpha,i)\sim (\alpha+1,i+c)$ if $i$ is even and 
$(\alpha,i)\sim (\alpha-1,i-c)$ if $i$ is odd. 
\end{definition}

\begin{figure}[t]
    \centering
    \includegraphics[width=15cm]{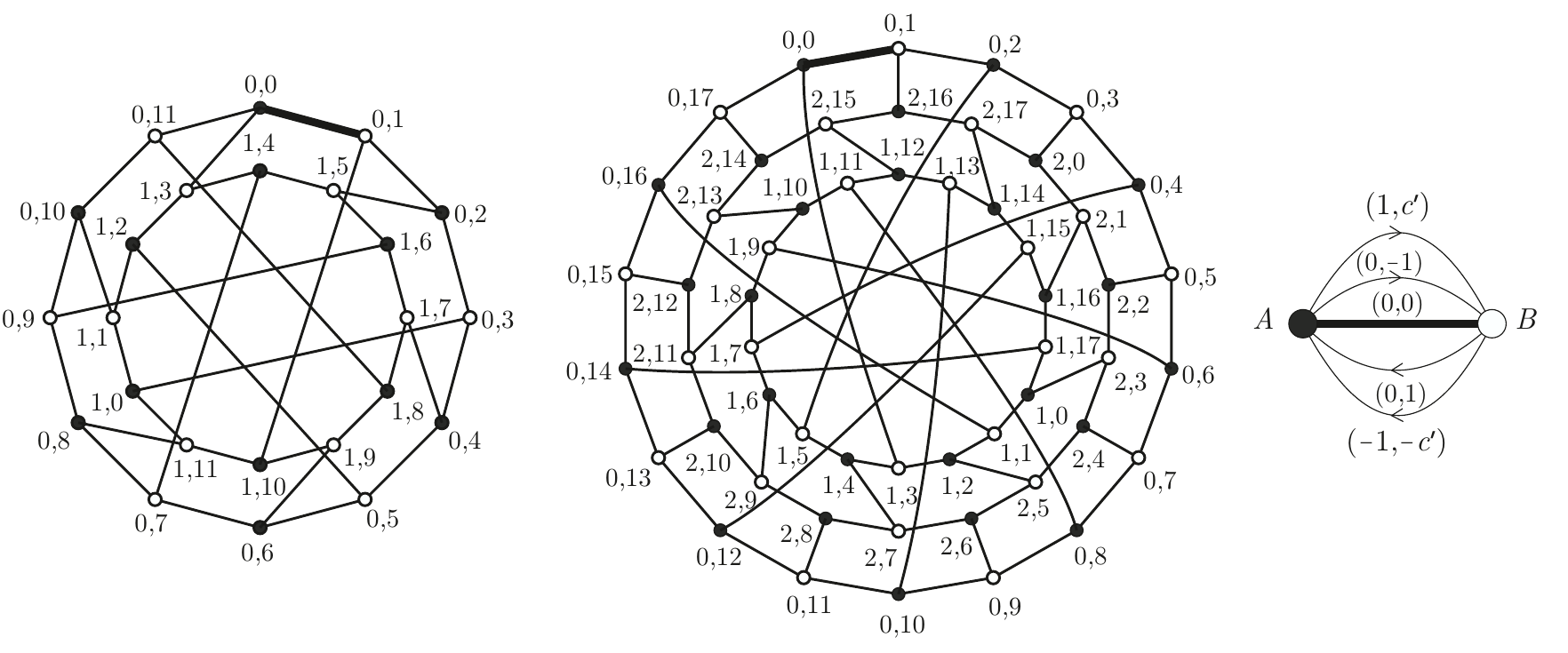}
    \caption{Left and middle: The chordal multi-ring graphs $CMR(2,12,3)$ and $CMR(3,18,3)$, where vertices $(x,y)$ are represented as $x,y$. Right: The base graph of the chordal multi-ring graph $CMR(m,n,c)$ on the Abelian group $\Z_{m}\times\Z_n$, where $c'=(c-1)/2$. The thick lines represent the edges with voltage $(0,0)$ in the base graph and in the first copies of the lift graphs.}
    \label{fig:cmr}
\end{figure}

Then, the graph $CMR(m,n,c)$ on $N=mn$ vertices is 3-regular and bipartite and consists of $m$ cycles of even length $n=2\nu$, together with some edges joining them. In particular, when $m=1$,  $CMR(1,n,c)\cong CR(n,c)$. For example, the chordal multi-ring graphs $CMR(2,12,3)$ and $CMR(3,18,3)$ are represented in Figure \ref{fig:cmr} (left and middle). 
As we show later, the first one is of special interest because all its eigenvalues are integers.

The adjacencies of the chordal $m$-ring
graphs follow the same planar pattern for every $m\ge 1$, see Figure \ref{fig:planar-pattern}. Then, their maximum numbers of vertices are again those in \eqref{boundCR}. As we show in what follows, the advantage of using more than one cycle is that, for an even diameter $k$, the maximum number of vertices can be attained. 

In Table \ref{tab:tabla-cmr}, we show the minimum diameter $k$ and chord $c$ for each number of vertices $N\le 138$ of a chordal multi-ring graph.
Table \ref{tab:tabla-CRM-repetida} provides the same results, but now, with only one value of $mn$ up to $N\le 354$.
The cases in which we get the maximum number of vertices for a given diameter are in boldface. As commented above, observe that, for an even diameter, the number of vertices of the chordal multi-ring graphs attains the maximum possible value. This can be proved in general.

\begin{proposition}
For an even diameter $k\ge 2$, the chordal $m$-ring graph $CMR(m,n,c)$ with $m=k/2$, $n=3k$, and $c=3$
has the maximum possible order $N=\frac{3}{2}k^2$.
\end{proposition}

\begin{figure}[t]
    \centering
    \includegraphics[width=15cm]{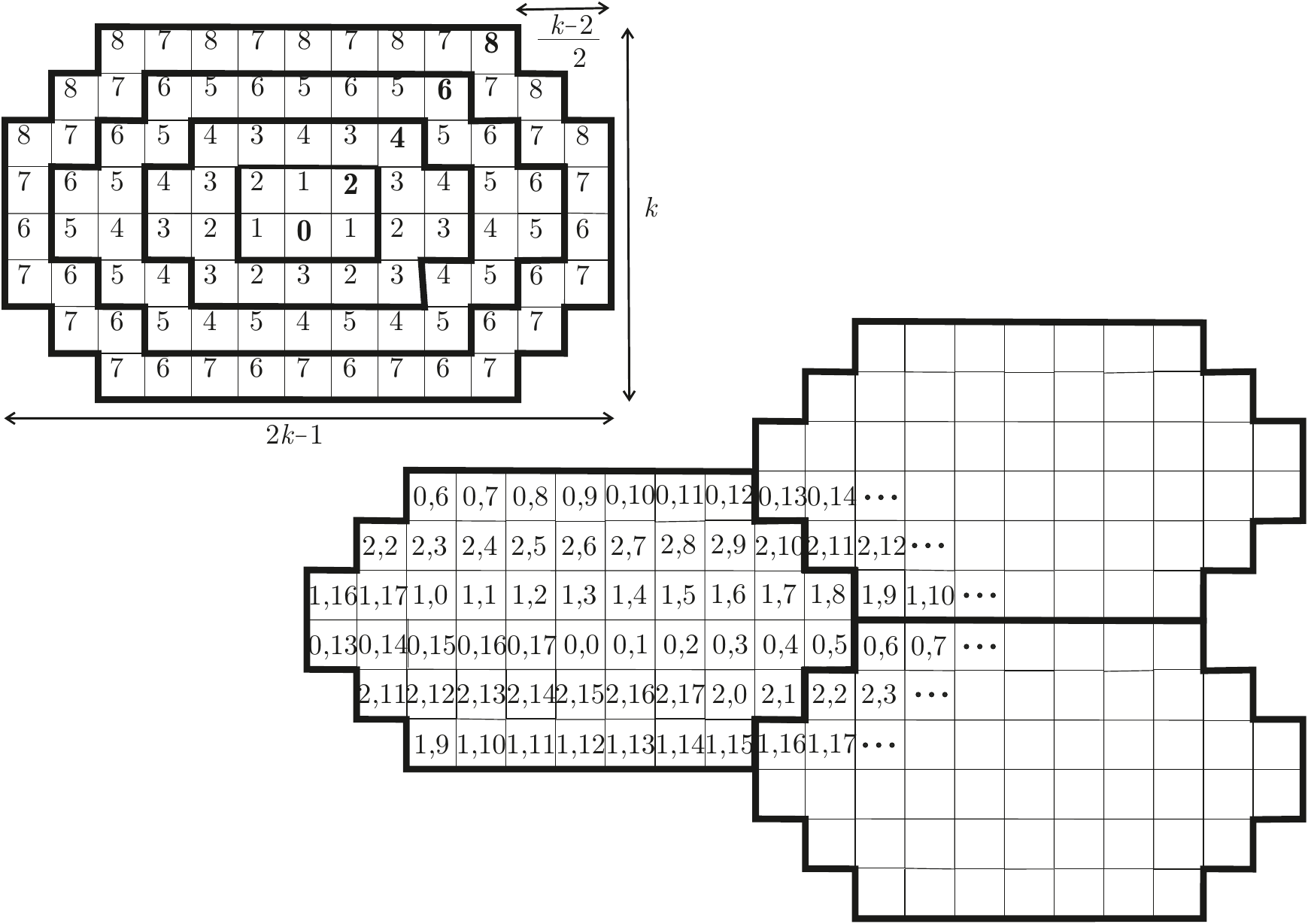}
    \caption{The optimal tiles for chordal multi-ring graphs with an even diameter $2,4,6,8$ (up).
    The planar pattern of $CMR(3,18,3)$ on the Abelian group $\Z_{3}\times\Z_{18}$ (down).}
    \label{fig:cmr-pattern}
\end{figure}

\begin{proof}
For every even $k$, an optimal tile with an area $3k^2/2$ periodically
tessellates the plane and corresponds to a chordal $m$-ring graph with $m = k/2$. 
See Figure \ref{fig:cmr-pattern} (up) for $k=2,4,6,8$, where the distances from vertex 0 are indicated. Then, from the distribution of the 0's, it can be checked that the corresponding graph has chord $c=3$. See Figure \ref{fig:cmr-pattern} 
 (down) for $k=6$ with the vertices labeled like in Definition \ref{def:CMR}. The obtained chordal $3$-ring graph is shown in Figure \ref{fig:cmr} (middle).
\end{proof}

Chordal multi-ring graphs can also be represented as lift graphs; see their base graph in Figure \ref{fig:cmr} (right). 
Then, we have the following result, which gives the eigenvalues of a chordal
$m$-ring graph $CMR(m, n, c)$.

\begin{proposition}
Given integers $m$, $n=2\nu$ (even),
and $c$ (odd), the eigenvalues of the chordal $m$-ring graph $CMR(m,n,c)$ are
\begin{align}
\lambda(r,s)_{1,2}&=\nonumber
\\ &\pm\sqrt{4\cos^2\left(\frac{s\pi}{\nu}\right)
+4\cos^2\left(\frac{r\pi}{m}+\frac{(c-1)s\pi}{n}\right)
+4\cos^2\left(\frac{r\pi}{m}+\frac{(c+1)s\pi}{n}\right)-3} 
\label{specCR(N,c)2}
\end{align}
for $r=0,1,\ldots,m-1$ and $s=0,1,\ldots,\nu-1$.
\end{proposition}

\begin{proof}
The base graph of $CMR(m,n,c)$ is shown in Figure \ref{fig:cmr} (right).
Then, its $(y,z)$-polynomial matrix is
$$
		\B(y,z)=
		\left(
		\begin{array}{cc}
		0 & 1+\frac{1}{z}+yz^{c'}  \\
		1+z+\frac{1}{yz^{c'}} & 0 
		\end{array}
		\right),
$$
where $c'=\frac{c-1}{2}$.
Then, the eigenvalues of $CMR(m,n,c)$ can be obtained as the eigenvalues of $\B(y,z)$ for every $y=e^{\frac{i2\pi}{m}r}$ with $\nu=n/2$, $r=0,\ldots,m-1$, and $z=e^{\frac{i2\pi}{\nu}s}$ for $s=0,\ldots,\nu-1$.
With  $\tau=1+\frac{1}{z}+yz^{c'}$, such eigenvalues are 
\begin{align*}
&\lambda(r)_{1,2}=\pm\sqrt{\tau\overline{\tau}}=\pm|\tau|\\
 &=\pm\sqrt{\left[1+\cos\left(\frac{s2\pi}{\nu}\right)+\cos\left(\frac{r2\pi}{m}+\frac{c's2\pi}{\nu}\right)\right]^2
 +\left[-\sin\left(\frac{s2\pi}{\nu}\right)+\sin\left(\frac{r2\pi}{m}+\frac{c's2\pi}{\nu}\right)\right]^2}
\end{align*}
and, operating with $c'=\frac{c-1}{2}$, we obtain \eqref{specCR(N,c)2}.
\end{proof}
Figure \ref{fig:cmr} shows the chordal multi-ring graphs $CMR(2, 12, 3)$ and $CMR(3, 18, 3)$.
Their eigenvalues are shown in Table  
\ref{tab:ev-cmr} with $y=e^{ri\frac{2\pi}{2}}$ and $z=e^{si\frac{2\pi}{6}}$. 
Notice that the graph $CMR(2,12,3)$ has an integral spectrum; see, for example, Ahmadi et al. \cite{aabs09}. As commented in that paper, such graphs play an important role in quantum
networks supporting the so-called perfect state transfer.

  \begin{table}[ht]
  \caption{The eigenvalues of the chordal multi-ring graphs  $CMR(2,12,3)$ and $CMR(3,18,3)$, respectively.}
		\label{tab:ev-cmr}
  \small
		\begin{center}
  \setlength{\tabcolsep}{5pt}
  \begin{tabular}{|c||c|c|c|c|c|c|}
	\hline
$r\setminus s$ & 0 & 1 & 2 & 3 & 4 & 5 \\
	\hline\hline
0 & $\pm 3$ & $\pm 2$ & $\pm 0$ & $\pm 1$ & $\pm 0$ & $\pm 2$ \\ 
     \hline
1 & $\pm 1$ & $\pm 2$ & $\pm 2$ & $\pm 1$ & $\pm 2$ & $\pm 2$ \\
\hline
				\end{tabular}
    \vskip .5cm
	\begin{tabular}{|c||c|c|c|c|c|c|c|c|c|}
	\hline
$r\setminus s$ & 0 & 1 & 2 & 3 & 4 & 5 & 6 & 7 & 8 \\
	\hline\hline
0 & $\pm 3$ & $\pm 0.8794$ & $\pm 2.532$ & $\pm 0$ & $\pm 1.348$ & $\pm 1.348$ & $\pm 0$ & $\pm 2.532$ & $\pm 0.8794$ \\ 
     \hline
1 & $\pm 1.732$ & $\pm 2.532$ & $\pm 1.348$ & $\pm 1.732$ & $\pm 0.8794$ & $\pm 0.8794$ & $\pm 1.732$ & $\pm 1.348$ & $\pm 2.532$ \\
     \hline
2 & $\pm 1.732$ & $\pm 2.532$ & $\pm 1.348$ & $\pm 1.732$ & $\pm 0.8794$ & $\pm 0.8794$ & $\pm 1.732$ & $\pm 1.348$ & $\pm 2.532$ \\
					\hline
				\end{tabular}
		\end{center}
		\end{table}


\section{Chordal ring mixed graphs} 
\label{sec:CRM}

In a mixed graph, there are edges (without direction) and arcs (with direction). See some recent results of some mixed graphs in 
Dalf\'o et al. \cite{deeft24} and 
Dalf\'o et al. \cite{deeflmt24}.

\begin{figure}[t]
    \centering
    \includegraphics[width=14cm]{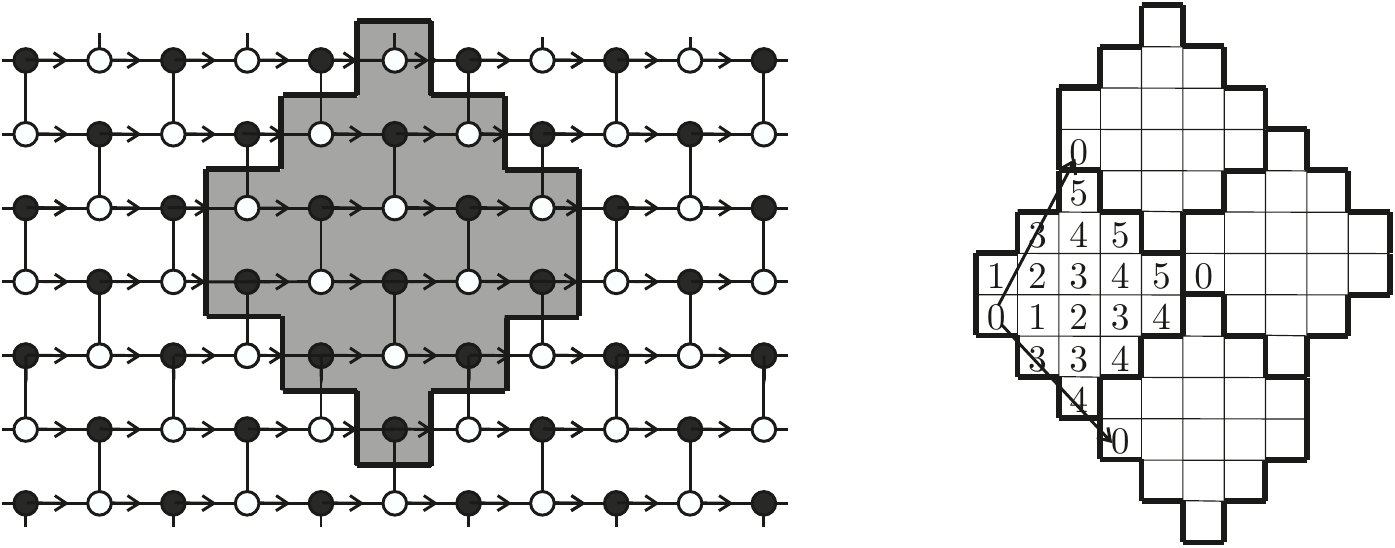}
    \caption{The optimal tiles for the chordal ring mixed graphs with a diameter $k=5$. The translation vectors are $(2,4)^{\top}$ and 
    $(3,-3)^{\top}$.}
    \label{fig:diameter5-mixed}
\end{figure}

\begin{definition}
\label{defCRM}
Let $N\ge 2$ and $c(<N)$  be, respectively, even and odd numbers. The chordal ring mixed graph $CRM(N,c)$ is a mixed graph with vertex set $V=\Z_N$ (all arithmetic is modulo $N$), with arcs $i\rightarrow i+1$ (forming a directed cycle) and edges $i\sim i+c$ if $i$ is even  (and, hence, $i\sim i-c$ if $i$ is odd, forming the `chords').
\end{definition}

See the example of the chordal ring mixed graph $CRM(10,3)$ in Figure \ref{fig:qCRM} (left).

\begin{figure}[!ht]
    \centering
    \includegraphics[width=14cm]{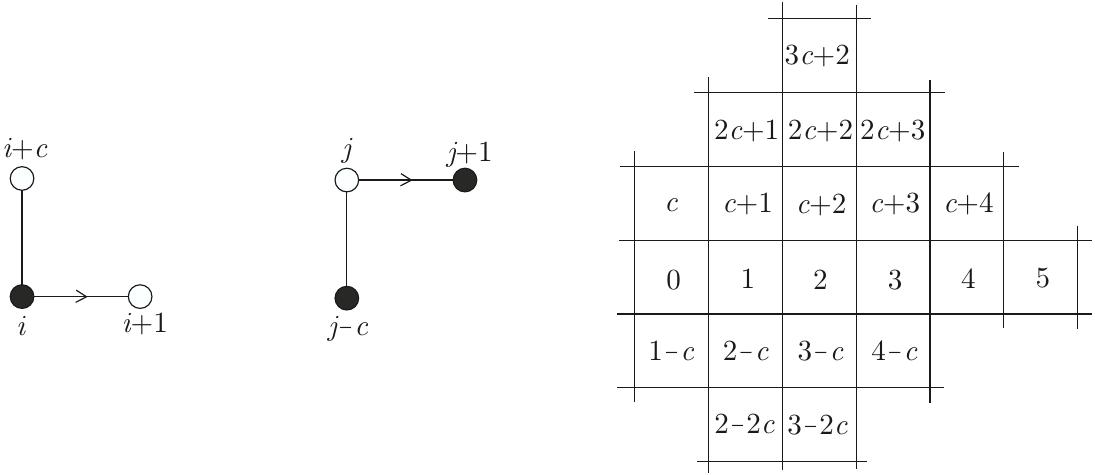}
    \caption{Adjacencies (for $i$ even and $j$ odd) and planar pattern of the vertices in the chordal ring mixed graph $CRM(N,c)$ (vertices at a distance at most five from vertex $0$).}
    \label{fig:planar-pattern-mixed}
\end{figure}

As in the case of chordal ring graphs $CR(N,c)$, if each vertex of $CRM(N,c)$ is represented by a numbered unit square modulo $N$, the vertices reached at a distance $0,1,2,\ldots$ from any given vertex can be arranged in a planar pattern, as shown in Figure \ref{fig:planar-pattern-mixed} (starting from vertex $0$).

Since there are at most $\ell+1$ vertices at a distance $\ell(>0)$ from vertex $0$, and the graph is bipartite,  the maximum number $N$ of vertices of a chordal ring mixed graph with a diameter $k$ is 
\begin{equation}
N(k)=\left\{ 
\begin{array}{ll}
\frac{1}{2}(k+1)^2 & \mbox{for $k$ odd,}\\[.1cm]
\frac{1}{2}k(k+2) & \mbox{for $k$ even.}
\end{array}
\right.
\label{boundCRM}
\end{equation}
Moreover, 
Dalf\'o et al. \cite{deeft24} showed that such a maximum can be attained when $k$ is odd but cannot be attained when $k(>2)$ is even.
We can raise the following conjecture from the reasoning in \cite{deeft24} and computer exploration shown in Table \ref{tab:tabla-CRM}.

\begin{conjecture}
\label{conj:CRM}
The maximum number $N$ of vertices of a chordal ring mixed graph $CRM(N,c)$ with an even diameter $k>2$ is $N=\frac{1}{2}k^2+2$ if $k\equiv 0 \mod 4$ (with $c=\frac{1}{4}n^2-\frac{1}{2}n+1$), and $N=k(\frac{k}{2}-1)+4$ if $k\equiv 2 \mod 4$. 
\end{conjecture}

As in the case of chordal ring graphs, the method used in 
\cite{deeft24} to find the optimal ring mixed graphs consists of the same steps (1) and  (2) from the Introduction. For example, Figure \ref{fig:diameter5-mixed} (right) shows the optimal tile for a diameter  $k=5$ and its periodic tessellation.

In Table \ref{tab:tabla-CRM}, we show the minimum diameter $k$ and chord $c$ for each number of vertices $N\le 526$ of a chordal ring mixed graph.
The cases in which we get the maximum number of vertices for a given diameter are in boldface.


\subsection{Chordal ring mixed graphs as lifts}

\begin{figure}[ht]
    \centering
    \includegraphics[width=12cm]{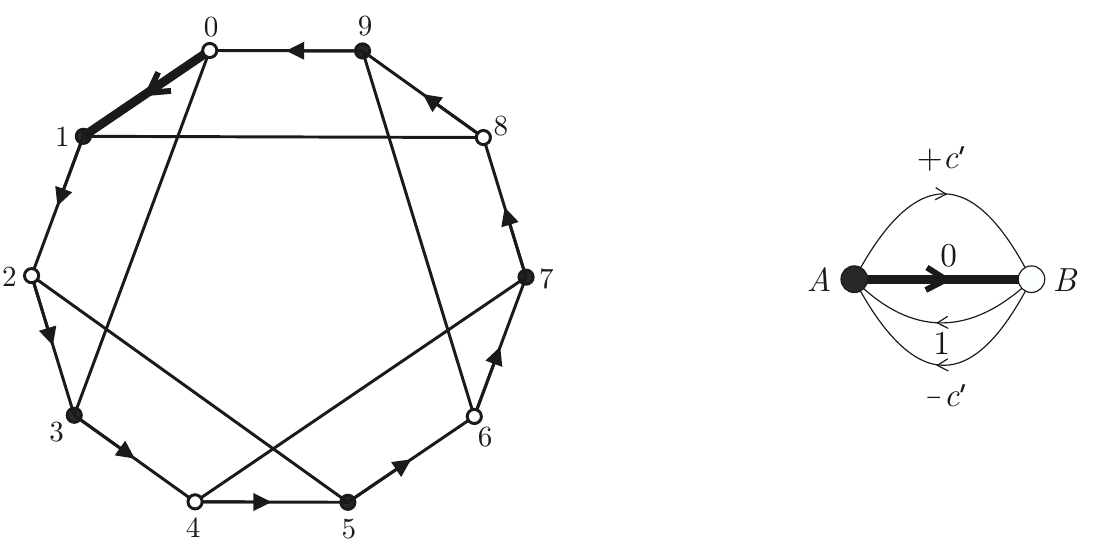}
    \caption{Left: The chordal ring mixed graph $CRM(10,3)$, 
    where vertices $(x,y)$ are denoted $x,y$. Right: The base graph of the chordal ring mixed graph $CRM(N,c)$ on the group $\Z_{N/2}$, where $c'=(c-1)/2$. The thick lines represent the arcs with voltage 0 in the base graph and in the first copy of the lift graph.}
    \label{fig:qCRM}
\end{figure}

The chordal ring mixed graph $CRM(N,c)$ can be seen as a lift of a base mixed graph on the group $\Z_{N/2}$, which is represented in Figure \ref{fig:qCRM} (right).

\begin{proposition}
Given integers $N=2n$ and $c$ (odd), the eigenvalues of the chordal ring graph mixed graphs $CRM(N,c)$ are
\begin{equation}
    \pm\lambda_{1,2}(r)=\sqrt{z^{c'}(z^{1+2c'}+z^{1+c'}+z^{c'}+1})\, z^{-c'},
    \label{specCRM(N,c)}
\end{equation}
where $z=e^{\frac{i2\pi}{n}r}$, $c'=\frac{c-1}{2}$, 
and $r=0,1,\ldots,n-1$.
\end{proposition}

\begin{proof}
First, note that $CRM(N,c)$ can be obtained as the lift of the base graph on the right side of Figure \ref{fig:qCRM}.
The polynomial matrix $\B(z)$, with $z=e^{\frac{i2\pi}{n}r}$, of such a base graph is 
\begin{equation*}
\B(z)=
\left(
\begin{array}{cc}
0 & 1+z^{c'}  \\
z+\frac{1}{z^{c'}} & 0 
\end{array}
\right),
\label{B(z)-CMR}
\end{equation*}
where $z=e^{\frac{i2\pi}{n}r}$ and $c'=(c-1)/2$. Then, the eigenvalues of  the lift $CRM(n,c)$ are the eigenvalues of $\B(z)$ for $r=0,1,\ldots,n-1$, given in \eqref{specCRM(N,c)}.
\end{proof}
Note that, in general, the matrix $\B(z)$ is not Hermitian; hence, the obtained eigenvalues are complex numbers.

\begin{figure}[t]
    \centering
    \includegraphics[width=6cm]{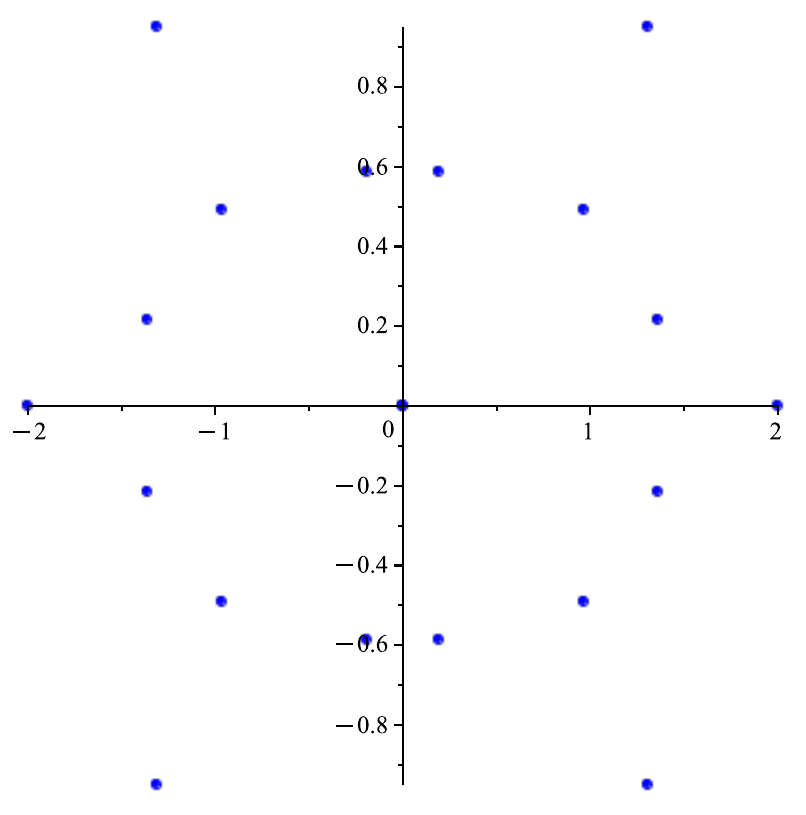}
    \caption{The complex eigenvalues of the chordal ring mixed graph $CRM(20,5)$.}
    \label{fig:ev-CRM}
\end{figure}

\begin{example}
Consider the case of $CRM(20,5)$, that is, the chordal ring mixed graph with 20 vertices and, for each vertex $i$, there are edges $i\sim i\pm 5$ and one arc $i\rightarrow i+1$  
(see Definition \ref{defCRM}). We obtain the eigenvalues in Table \ref{tab:ev-CRM}, which are represented in Figure \ref{fig:ev-CRM}. 

  \begin{table}[ht]
  \caption{All the eigenvalues of the chordal ring mixed graph $CRM(20,5)$.}
		\label{tab:ev-CRM}
		\begin{center}
	\begin{tabular}{|c|c|}
	\hline
	\rule{0pt}{3ex} $\zeta=e^{i\frac{2\pi}{10}}$, $z=\zeta^r$  & $\lambda$ \\[0.5ex]
	\hline\hline
	$\spec(\B(1))$ &  $\pm 2$ \\
     \hline
    $\spec(\B(\zeta))$ &  $\pm(1.362+0.2158i)$\\
     \hline
     $\spec(\B(\zeta^2))$ &  $\pm (0.1913-0.5879i)$ \\
      \hline
    $\spec(\B(\zeta^3))$ &  $\pm(0.9661+0.4922i)$ \\
      \hline
      $\spec(\B(\zeta^4))$ &  $\pm (1.309+0.9511i)$ \\
      \hline
      $\spec(\B(\zeta^5))$ &  $\pm 0$ \\
      \hline
      $\spec(\B(\zeta^6))$ &  $\pm (1.309-0.9511i)$ \\
      \hline
      $\spec(\B(\zeta^7))$ &  $\pm (0.9661-0.4922i)$ \\
       \hline
      $\spec(\B(\zeta^8))$ &  $\pm (0.1913+0.5879i)$ \\
     \hline
     $\spec(\B(\zeta^9))$ &   $\pm(1.362-0.2158i)$\\
					\hline
				\end{tabular}
		\end{center}
		\end{table}
	\end{example}

\section{Chordal multi-ring mixed graphs} 
\label{sec:CMRM}

As in the case of chordal ring graphs, we can consider the mixed version of the chordal multi-ring graphs, denoted $CMRM(m,n,c)$, which, as mixed graphs, have edges and arcs.

\begin{definition}
Given positive integers $m$, $n$
(even) and $c(>1)$ (odd), the chordal 
$m$-ring mixed graph $CMRM(m,n,c)$ has vertices labeled with the elements of the Abelian group $\Z_m\times \Z_n$, arcs $(\alpha,i)\rightarrow (\alpha,i+1)$ for every $\alpha\in \Z_m$ and $i\in \Z_n$, and edges $(\alpha,i)\sim (\alpha+1,i+c)$ if $i$ is even and 
$(\alpha,i)\sim (\alpha-1,i-c)$ if $i$ odd. 
\end{definition}

\begin{figure}[t]
    \centering
    \includegraphics[width=12cm]{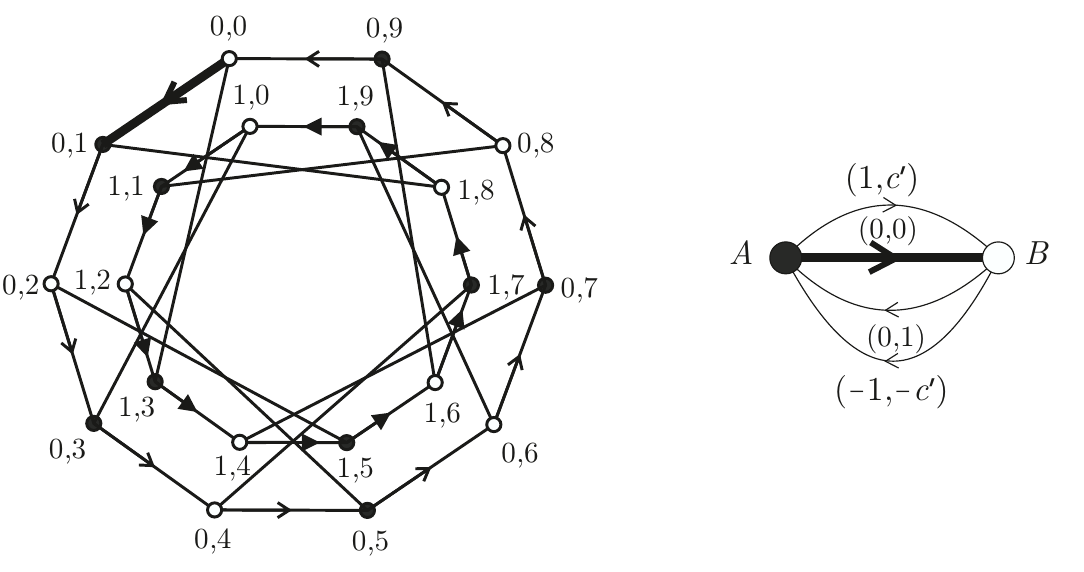}
    \caption{The chordal multi-ring mixed graph $CMRM(2,10,3)$ (left) and the base graph of the chordal multi-ring mixed graph $CMRM(m,n,c)$ (right). The thick lines represent the arcs with voltage $(0,0)$ in the base graph and in the first copy of the lift graph.}
  \label{fig:base-CMRM}
\end{figure}

See the example of the chordal multi-ring mixed graph $CMRM(2,10,3)$ in Figure \ref{fig:base-CMRM} (on the left), and its eigenvalues in Table \ref{tab:ev-CMRM}.

Similar to the case of chordal multi-ring graphs, the adjacencies of the chordal $m$-ring mixed 
graphs follow the same planar pattern for every $m\ge 1$, see Figure \ref{fig:planar-pattern-mixed}. Then, their maximum numbers of vertices are again those in \eqref{boundCRM}. Now, by using more than one cycle, we can improve the number of vertices reached for an even diameter $k$, but, in this case, without getting the possible maximum value in \eqref{boundCRM}. More precisely, as shown in Conjecture \ref{conj:CMRM} below, we get $N=\frac{1}{2}k^2+2$ when $k\equiv 2 \mod 4$ (instead of $N=k(\frac{k}{2}-1)+4$
in Conjecture \ref{conj:CRM}).

In Table \ref{tab:tabla-cmrm}, we show the minimum diameter $k$ and chord $c$ for each number of vertices $N\le 354$ of a chordal multi-ring mixed graph.
The cases in which we get the maximum number of vertices for a given diameter are in boldface.
From these values, we pose the following conjecture.
 
\begin{conjecture}
\label{conj:CMRM}
The maximum number $N$ of vertices of a chordal $m$-ring mixed graph with an even diameter $k\equiv 2 \mod 4$ is $N=\frac{1}{2}k^2+2$,
and this value is attained with $m=2$ cycles and chord $c=k/2$.
\end{conjecture}

We have the following result concerning the spectrum of the chordal multi-ring mixed graphs.



\begin{proposition}
Given integers $m$, $n$ (even), and $c$ (odd), the eigenvalues of the chordal multi-ring mixed graph $CMRM(m,n,c)$ are
\begin{equation}
\pm\lambda_{1,2}(r)=\sqrt{yz^{c'}(y^2 z^{1+2c'}+yz^{1+c'}+yz^{c'}+1)}\,\, y^{-1}z^{-c'},
\label{specCMRM(N,c)}
\end{equation}
where $y=e^{\frac{i2\pi}{m}r}$, for $r=0,1,\ldots,m-1$, $z=e^{\frac{i2\pi}{n}s}$, $s=0,1,\ldots,n-1$, and $c'=\frac{c-1}{2}$.
\end{proposition}

\begin{proof}
In this case, $CMRM(m,n,c)$ can be obtained as the lift of the base graph on the right side of Figure \ref{fig:base-CMRM}.
The polynomial matrix $\B(y,z)$ of such a base graph is 
\begin{equation*}
\B(y,z)=
\left(
\begin{array}{cc}
0 & 1+yz^{c'}  \\
z+\frac{1}{yz^{c'}} & 0 
\end{array}
\right),
\end{equation*}
with $y=e^{\frac{i2\pi}{m}r}$ and $z=e^{\frac{i2\pi}{n}r}$. Then, the eigenvalues of  the lift $CMRM(m,n,c)$ are the eigenvalues of $\B(y,z)$ for $r=0,1,\ldots,n-1$, and $s=0,1,\ldots,n-1$ given in \eqref{specCMRM(N,c)}.
\end{proof}

  \begin{table}[ht]
  \caption{All the eigenvalues of the chordal multi-ring mixed graph $CMRM(2,10,3)$.}
		\label{tab:ev-CMRM}
		\begin{center}
\begin{tabular}{|c||c|c|}
	\hline
$s\setminus r$ & 0 & 1 \\
	\hline\hline
0 & $\pm 2$ & $\pm 0$ \\ 
     \hline
1 &  $\pm (0.9511+0.3090i)$ & $\pm (0.8419+i0.1337)$  \\
\hline
2 &  $\pm (1.733+0.2744i)$ & $\pm (1.422+i0.4620)$  \\
\hline
3 &  $\pm (0.3870-0.75941i)$ & $\pm (1.563+i0.7965)$  \\
\hline
4 &  $\pm (-0.5878+0.8090i)$ & $\pm(1.210+i0.8790)$  \\
					\hline
				\end{tabular}
    \end{center}
		\end{table}

\section{Discussion}

In this paper, we have obtained the following results:
\begin{enumerate}
\item We have introduced the families of chordal multi-ring ($CMR$), chordal ring mixed 
($CMR$), and chordal multi-ring mixed ($CMRM$) graphs from the well-known chordal ring 
($CR$) graphs.
\item For these four families, we have obtained the maximum number of vertices
for every given diameter by using the approach of plane tessellations and their corresponding lattices. Moreover, we computed the minimum diameter for any value of the number of vertices. Actually, for our families of graphs, these are optimization results that solve the degree/diameter problem and the degree/number of vertices problem, which are very well-known problems in the literature. See, for example, the comprehensive survey by Miller and \v{S}ir\'a\v{n} \cite{ms13}.
\item Finally, we have obtained the spectra of these four families.
Here, we used the theory of lifts of some base graphs on Abelian groups to derive closed formulas for the eigenvalues of our families of graphs.
\end{enumerate}

\section*{Statements and Declarations}

The authors declare no conflicts of interest.\\
The funders had no role in the design of the study, in the collection, analyses, or interpretation of data, in the writing of the manuscript, or in the decision to publish the results.
All authors have contributed equally in conceptualization, methodology, software, validation, formal analysis, investigation, resources, data curation, writing, and reviewing.\\
All authors have read and agreed to the published version of the manuscript.\\
All data generated with this paper is included. 


\section*{Appendix A}

In the following table, we give a list of the main symbols used, together with their meaning.

\begin{table}[!ht]
\caption{Nomenclature}
\begin{center}
\begin{tabular}{|c|l|}
\hline
Symbol & Description\\
\hline\hline 
$\alpha$ & Voltage assignment\\
$\A$ & Adjacency matrix\\
$\B(z)$ & Polynomial matrix with complex entries\\
$c$ & Chord\\
$CR$ & Chordal ring graph\\
$CMR$ & Chordal multi-ring graph\\
$CRM$ & Chordal ring mixed graph \\
$CMRM$ & Chordal multi-ring mixed graph\\
$E$ & Set of edges\\
$\Gamma$ & Graph \\
$\Gamma^{\alpha}$ & Lift graph of $\Gamma$\\
$G$ & Group\\
$k$ & Diameter of the graph\\
$\lambda(r)$  & Eigenvalue of the matrix $\B(\zeta^r)$\\
$\M$ & Integral matrix\\
$N$ & Number of vertices\\
$\mathbb{R}_k[z]$ & Ring of complex polynomials with degree $\le k$\\
$V$ & Set of vertices \\
$\zeta$ & Primitive $n$-th root of unity\\
$\Z_N$ & Cyclic group of integers modulo $N$\\
$\Z^2\M$ & Lattice generated by the matrix $\M$\\
\hline
\end{tabular}
\end{center}
\end{table}

\section*{Appendix B}
This section shows the tables of our different families of chordal ring graphs and mixed graphs.

They have been obtained by an exhaustive computer search, and their principal values agree with our theoretical results and conjectures.

  \begin{table}[t]
  \caption{The minimum diameter $k$ for each number of vertices $N\le 528$ in the chordal ring graphs $CR(N,c)$ with (minimum) chord $c$. In boldface, there are the cases with a maximum number of vertices for a given diameter.}
		\label{tab:tabla-cr}
  \small
		\begin{center}
\setlength{\tabcolsep}{3pt}
\begin{tabular}{|c|c|c||c|c|c||c|c|c||c|c|c||c|c|c||c|c|c|}
\hline
$N$ & $c$ & $k$ & $N$ & $c$ & $k$ & $N$ & $c$ & $k$ & $N$ & $c$ & $k$ & $N$ & $c$ & $k$ & $N$ & $c$ & $k$\\
\hline
  {\bf 2} &   {\bf 1} &   {\bf 1} &	 90 &  11 &   9 &	178 &  49 &  12 &	266 &  41 &  14 &	354 &  47 &  16 &	442 &  47 &  19 \\
  4 &   1 &   2 &	 92 &  11 &   9 &	180 &  33 &  12 &	268 &  41 &  14 &	356 &  47 &  17 &	444 &  47 &  19 \\
  {\bf 6} &   {\bf 3} &  {\bf 2} &	 94 &  11 &   9 &	{\bf 182} &  {\bf 33} &  {\bf 11} &	270 &  21 &  15 &	358 &  55 &  17 &	446 &  47 &  19 \\
  8 &   3 &   3 &	 96 &  11 &   9 &	184 &  33 &  12 &	272 &  23 &  15 &	360 &  67 &  17 &	448 &  47 &  19 \\
 10 &   3 &   3 &	 98 &  11 &   9 &	186 &  51 &  12 &	274 &  37 &  15 &	362 &  43 &  17 &	450 &  53 &  18 \\
 12 &   5 &   3 &	100 &  13 &   9 &	188 &  15 &  13 &	276 &  37 &  15 &	364 &  43 &  17 &	452 &  53 &  18 \\
 {\bf 14} &   {\bf 5} &  {\bf 3} &	102 &  13 &   9 &	190 &  79 &  12 &	278 &  37 &  15 &	366 &  43 &  17 &	454 &  53 &  19 \\
 16 &   5 &   4 &	104 &  23 &   9 &	192 &  35 &  12 &	{\bf 280} &  {\bf 43} &  {\bf 14} &	{\bf 368} &  {\bf 49} &  {\bf 16} &	456 &  81 &  19 \\
 18 &   5 &   4 &	106 &  23 &   9 &	194 &  35 &  12 &	282 &  37 &  15 &	370 &  43 &  17 &	458 &  61 &  19 \\
 {\bf 20} &   {\bf 5} &   {\bf 4} &	108 &  41 &   9 &	196 &  17 &  13 &	284 &  61 &  15 &	372 & 163 &  17 &	460 &  97 &  19 \\
 22 &   5 &   5 &	110 &  11 &  10 &	198 &  17 &  13 &	286 &  21 &  16 &	374 &  27 &  18 &	462 &  49 &  19 \\
 24 &   5 &   5 &	112 &  25 &   9 &	200 &  17 &  13 &	288 &  77 &  15 &	376 &  81 &  17 &	464 &  49 &  19 \\
 26 &   5 &   5 &	114 &  25 &   9 &	202 &  17 &  13 &	290 &  39 &  15 &	378 &  67 &  17 &	466 &  49 &  19 \\
 28 &   7 &   5 &	116 &  13 &  10 &	{\bf 204} &  {\bf 37} &  {\bf 12} &	292 &  39 &  15 &	380 &  45 &  17 &	{\bf 468} &  {\bf 55} &  {\bf 18} \\
 30 &   7 &   5 &	118 &  15 &  10 &	206 &  19 &  13 &	294 &  39 &  15 &	382 &  45 &  17 &	470 &  49 &  19 \\
 32 &   7 &   5 &	120 &  27 &  10 &	208 &  57 &  13 &	296 &  39 &  15 &	384 &  45 &  17 &	472 & 131 &  19 \\
 34 &   7 &   5 &	{\bf 122} &  {\bf 27} &   {\bf 9} &	210 &  77 &  13 &	298 &  39 &  16 &	386 &  45 &  17 &	474 &  63 &  20 \\
 36 &   7 &   6 &	124 &  27 &  10 &	212 &  33 &  13 &	300 &  69 &  15 &	388 &  45 &  18 &	476 &  75 &  19 \\
 {\bf 38} &  {\bf 15} &   {\bf 5} &	126 &  13 &  11 &	214 &  33 &  13 &	302 &  65 &  15 &	390 &  69 &  17 &	478 & 141 &  19 \\
 40 &   7 &   6 &	128 &  13 &  11 &	216 &  33 &  13 &	304 &  41 &  16 &	392 &  69 &  18 &	480 &  85 &  19 \\
 42 &   7 &   6 &	130 &  29 &  10 &	218 &  33 &  13 &	306 &  41 &  15 &	394 &  85 &  18 &	482 &  51 &  19 \\
 44 &   9 &   6 &	132 &  29 &  10 &	220 &  17 &  14 &	308 &  41 &  15 &	396 &  47 &  18 &	484 &  51 &  19 \\
 46 &   7 &   7 &	134 &  13 &  11 &	222 &  51 &  13 &	310 &  41 &  15 &	398 &  47 &  17 &	486 &  51 &  19 \\
 {\bf 48} &  {\bf 19} &   {\bf 6} &	136 &  13 &  11 &	224 &  19 &  14 &	312 &  67 &  15 &	400 &  47 &  17 &	488 &  51 &  19 \\
 50 &   7 &   7 &	138 &  13 &  11 &	226 &  35 &  13 &	314 &  59 &  16 &	402 &  47 &  17 &	490 &  51 &  20 \\
 52 &   9 &   7 &	{\bf 140} &  {\bf 31} &  {\bf 10} &	228 &  35 &  13 &	316 &  59 &  16 &	404 & 119 &  17 &	492 &  87 &  19 \\
 54 &   9 &   7 &	142 &  13 &  11 &	230 &  35 &  13 &	318 &  49 &  16 &	406 & 113 &  18 &	494 & 137 &  19 \\
 56 &   9 &   7 &	144 &  15 &  11 &	232 &  89 &  13 &	320 &  43 &  16 &	408 & 109 &  18 &	496 & 183 &  20 \\
 58 &   9 &   7 &	146 &  15 &  11 &	234 &  17 &  15 &	322 &  43 &  15 &	410 & 121 &  18 &	498 &  75 &  20 \\
 60 &   9 &   7 &	148 &  15 &  11 &	236 &  43 &  14 &	324 &  43 &  15 &	412 &  55 &  18 &	500 &  53 &  20 \\
 62 &   9 &   7 &	150 &  17 &  11 &	238 &  37 &  14 &	326 &  43 &  16 &	414 &  49 &  18 &	502 &  53 &  19 \\
 64 &  11 &   7 &	152 &  27 &  11 &	240 &  37 &  13 &	328 &  21 &  17 &	416 &  49 &  17 &	504 &  53 &  19 \\
 66 &  19 &   7 &	154 &  13 &  12 &	242 &  37 &  13 &	330 &  71 &  16 &	418 &  49 &  17 &	506 &  53 &  19 \\
 68 &  19 &   7 &	156 &  43 &  11 &	244 &  37 &  14 &	332 &  51 &  16 &	420 &  49 &  18 &	508 & 141 &  19 \\
 70 &   9 &   8 &	158 &  29 &  11 &	246 &  53 &  14 &	334 &  21 &  17 &	422 &  79 &  18 &	510 &  27 &  21 \\
 72 &   9 &   8 &	160 &  29 &  11 &	248 &  45 &  14 &	336 &  45 &  16 &	424 &  75 &  18 &	512 &  77 &  20 \\
 {\bf 74} &  {\bf 21} &   {\bf 7} &	162 &  29 &  11 &	250 &  17 &  15 &	{\bf 338} &  {\bf 45} &  {\bf 15} &	426 &  27 &  19 &	514 &  77 &  20 \\
 76 &  11 &   8 &	164 &  45 &  11 &	252 &  39 &  14 &	340 &  45 &  16 &	428 &  57 &  18 &	516 & 143 &  20 \\
 78 &  29 &   8 &	166 &  17 &  12 &	{\bf 254} &  {\bf 39} &  {\bf 13} &	342 &  23 &  17 &	430 &  29 &  19 &	518 &  61 &  20 \\
 80 &  23 &   8 &	168 &  17 &  12 &	256 &  39 &  14 &	344 &  25 &  17 &	432 &  51 &  18 &	520 &  55 &  20 \\
 82 &  23 &   8 &	170 &  31 &  11 &	258 &  99 &  14 &	346 &  25 &  17 &	{\bf 434} &  {\bf 51} &  {\bf 17} &	522 &  55 &  19 \\
 84 &  11 &   9 &	172 &  31 &  11 &	260 &  19 &  15 &	348 &  41 &  17 &	436 &  51 &  18 &	524 &  55 &  19 \\
 86 &  11 &   9 &	174 &  31 &  12 &	262 &  19 &  15 &	350 & 153 &  16 &	438 & 129 &  18 &	526 &  55 &  20 \\
 {\bf 88} &  {\bf 25} &  {\bf 8} &	176 &  39 &  12 &	264 &  19 &  15 &	352 &  47 &  16 &	440 &  69 &  19 &	528 &  55 &  21 \\
\hline
		\end{tabular}
		\end{center}
		\end{table}

  \begin{table}[t]
  \caption{The minimum diameter $k$ in the chordal multi-ring graphs $CRM(m,n,c)$ with different values of $m$ and $n$. In boldface, there are the cases with a maximum number of vertices for a given diameter.}
		\label{tab:tabla-cmr}
  \small
		\begin{center}
\setlength{\tabcolsep}{3pt}
\begin{tabular}{|c|c|c|c|c||c|c|c|c|c||c|c|c|c|c||c|c|c|c|c|}
\hline
$N$ & $m$ & $n$ & $c$ & $k$ & $N$ & $m$ & $n$ & $c$ & $k$ & $N$ & $m$ & $n$ & $c$ & $k$ & $N$ & $m$ & $n$ & $c$ & $k$ \\
\hline
4 &   1 &   4&   3&   2 &46 &   1 &  46&   7&   7 &80 &   1 &  80&  23&   8 &108 &   9 &  12&   3&  10 \\
{\bf 6} &   {\bf 1} &   {\bf 6}&  {\bf 3}&   {\bf 2} &48 &   1 &  48&  19&   6 &80 &   2 &  40&   5&   8 &110 &   1 & 110&  11&  10 \\
8 &   1 &   8&   3&   3 &48 &   2 &  24&   5&   7 &80 &   4 &  20&   3&   8 &110 &   5 &  22&   3&  10 \\
8 &   2 &   4&   3&   3 &48 &   3 &  16&   3&   6 &80 &   5 &  16&   5&   9 &112 &   1 & 112&  25&   9 \\
10 &   1 &  10&   3&   3 &48 &   4 &  12&   5&   7 &80 &   8 &  10&   3&   9 &112 &   2 &  56&   7&  10 \\
12 &   1 &  12&   5&   3 &48 &   6 &   8&   3&   7 &82 &   1 &  82&  23&   8 &112 &   4 &  28&   3&  10 \\
12 &   2 &   6&   5&   3 &50 &   1 &  50&   7&   7 &84 &   1 &  84&  11&   9 &112 &   7 &  16&   3&  10 \\
12 &   3 &   4&   3&   3 &50 &   5 &  10&   3&   7 &84 &   2 &  42&   5&   9 &112 &   8 &  14&   3&  10 \\
{\bf 14} &   {\bf 1} &  {\bf 14}&   {\bf 5}&   {\bf 3} &52 &   1 &  52&   9&   7 &84 &   3 &  28&   3&   9 &114 &   1 & 114&  25&   9 \\
16 &   1 &  16&   5&   4 &52 &   2 &  26&   5&   7 &84 &   6 &  14&   3&   9 &114 &   3 &  38&  17&   9 \\
16 &   2 &   8&   3&   4 &54 &   1 &  54&   9&   7 &84 &   7 &  12&   3&   9 &116 &   1 & 116&  13&  10 \\
16 &   4 &   4&   3&   4 & {\bf 54} & {\bf 3} & {\bf 18} & {\bf 3} & {\bf 6} &86 &   1 &  86&  11&   9 &116 &   2 &  58&  11&  10 \\
18 &   1 &  18&   5&   4 &56 &   1 &  56&   9&   7 &88 &   1 &  88&  25&   8 &118 &   1 & 118&  15&  10 \\
18 &   3 &   6&   3&   4 &56 &   2 &  28&   5&   7 &88 &   2 &  44&   5&  10 &120 &   1 & 120&  27&  10 \\
20 &   1 &  20&   5&   4 &56 &   4 &  14&   5&   7 &88 &   4 &  22&   3&   8 &120 &   2 &  60&   7&  11 \\
20 &   2 &  10&   3&   4 &56 &   7 &   8&   3&   7 &90 &   1 &  90&  11&   9 &120 &   3 &  40&   5&  10 \\
22 &   1 &  22&   5&   5 &58 &   1 &  58&   9&   7 &90 &   3 &  30&   3&   9 &120 &   4 &  30&  11&  10 \\
24 &   1 &  24&   5&   5 &60 &   1 &  60&   9&   7 &90 &   5 &  18&   5&   9 &120 &   5 &  24&   3&  10 \\
{\bf 24} & {\bf 2} &  {\bf 12} & {\bf 3} & {\bf 4} &60 &   2 &  30&   5&   7 &90 &   9 &  10&   3&   9 &120 &   6 &  20&   5&  11 \\
24 &   3 &   8&   3&   5 &60 &   3 &  20&   3&   7 &92 &   1 &  92&  11&   9 &120 &  10 &  12&   3&  11 \\
24 &   4 &   6&   3&   5 &60 &   5 &  12&   3&   7 &92 &   2 &  46&   7&   9 &{\bf 122} & {\bf 1} & {\bf 122}& {\bf 27}& {\bf 9} \\
26 &   1 &  26&   5&   5 &60 &   6 &  10&   3&   7 &94 &   1 &  94&  11&   9 &124 &   1 & 124&  27&  10 \\
28 &   1 &  28&   7&   5 &62 &   1 &  62&   9&   7 &96 &   1 &  96&  11&   9 &124 &   2 &  62&  19&  10 \\
28 &   2 &  14&   3&   5 &64 &   1 &  64&  11&   7 &96 &   2 &  48&   7&   9 &126 &   1 & 126&  13&  11 \\
30 &   1 &  30&   7&   5 &64 &   2 &  32&   5&   8 &96 &   3 &  32&   7&   9 &126 &   3 &  42&   5&  10 \\
30 &   3 &  10&   5&   5 &64 &   4 &  16&   3&   8 &{\bf 96} & {\bf 4} & {\bf 24} & {\bf  3}& {\bf 8} &126 &   7 &  18&   3&  11 \\
30 &   5 &   6&   3&   5 &64 &   8 &   8&   3&   8 &96 &   6 &  16&   7&   9 &126 &   9 &  14&   3&  11 \\
32 &   1 &  32&   7&   5 &66 &   1 &  66&  19&   7 &96 &   8 &  12&   5&   9 &128 &   1 & 128&  13&  11 \\
32 &   2 &  16&   3&   6 &66 &   3 &  22&  11&   7 &98 &   1 &  98&  11&   9 &128 &   2 &  64&  13&  10 \\
32 &   4 &   8&   3&   5 &68 &   1 &  68&  19&   7 &98 &   7 &  14&   3&   9 &128 &   4 &  32&   3&  12 \\
34 &   1 &  34&   7&   5 &68 &   2 &  34&  11&   7 &100 &   1 & 100&  13&   9 &128 &   8 &  16&   3&  11 \\
36 &   1 &  36&   7&   6 &70 &   1 &  70&   9&   8 &100 &   2 &  50&   7&   9 &130 &   1 & 130&  29&  10 \\
36 &   2 &  18&   3&   6 &70 &   5 &  14&   5&   8 &100 &   5 &  20&   3&  10 &130 &   5 &  26&   3&  10 \\
36 &   3 &  12&   3&   6 &70 &   7 &  10&   3&   8 &100 &  10 &  10&   3&  10 &132 &   1 & 132&  29&  10 \\
36 &   6 &   6&   3&   6 &72 &   1 &  72&   9&   8 &102 &   1 & 102&  13&   9 &132 &   2 &  66&  13&  10 \\
{\bf 38} & {\bf 1} & {\bf 38} & {\bf 15} & {\bf 5} &72 &   2 &  36&   5&   8 &102 &   3 &  34&   7&   9 &132 &   3 &  44&   5&  10 \\
40 &   1 &  40&   7&   6 &72 &   3 &  24&   3&   9 &104 &   1 & 104&  23&   9 &132 &   6 &  22&   5&  11 \\
40 &   2 &  20&   3&   6 &72 &   4 &  18&   3&   8 &104 &   2 &  52&   7&   9 &132 &  11 &  12&   3&  11 \\
40 &   4 &  10&   3&   6 &72 &   6 &  12&   3&   8 &104 &   4 &  26&   3&   9 &134 &   1 & 134&  13&  11 \\
40 &   5 &   8&   3&   6 &{\bf 74} & {\bf 1} & {\bf 74}& {\bf 21}& {\bf 7} &106 &   1 & 106&  23&   9 &136 &   1 & 136&  13&  11 \\
42 &   1 &  42&   7&   6 &76 &   1 &  76&  11&   8 &108 &   1 & 108&  41&   9 &136 &   2 &  68&   7&  11 \\
42 &   3 &  14&   3&   6 &76 &   2 &  38&   5&   8 &108 &   2 &  54&   7&   9 &136 &   4 &  34&  13&  11 \\
44 &   1 &  44&   9&   6 &78 &   1 &  78&  29&   8 &108 &   3 &  36&   5&  10 &138 &   1 & 138&  13&  11 \\
44 &   2 &  22&   7&   6 &78 &   3 &  26&  13&   8 &108 &   6 &  18&   5&  10 &138 &   3 &  46&   5&  11 \\
\hline
\end{tabular}
		\end{center}
		\end{table}

  \begin{table}[t] 
  \caption{The minimum diameter $k$ in the chordal multi-ring graphs $CMR(m,n,c)$ with minimum values $m$ and $n$ as Table \ref{tab:tabla-cmr}, but now with only one value of $m$ and $n$.
In boldface, there are the cases with a maximum number of vertices for a given diameter.}
		\label{tab:tabla-CRM-repetida}
  \small
		\begin{center}
\setlength{\tabcolsep}{3pt}
\begin{tabular}{|c|c|c|c|c||c|c|c|c|c||c|c|c|c|c||c|c|c|c|c|}
\hline
$N$ & $m$ & $n$ & $c$ & $k$ & $N$ & $m$ & $n$ & $c$ & $k$ & $N$ & $m$ & $n$ & $c$ & $k$ & $N$ & $m$ & $n$ & $c$ & $k$ \\
\hline
4 &   1 &   4&   3&   2 &92 &   1 &  92&  11&   9 &180 &   1 & 180&  33&  12 &268 &   1 & 268&  41&  14 \\
{\bf 6} & {\bf 1} & {\bf 6} & {\bf 3} & {\bf 2} & 94 & 1 &  94&  11&   9 &{\bf 182} & {\bf 1} & {\bf 182} & {\bf 33}& {\bf 11} &270 &   1 & 270&  21&  15 \\
8 &   1 &   8&   3&   3 &{\bf 96} & {\bf 4} & {\bf 24}& {\bf 3}& {\bf 8} &184 &   1 & 184&  33&  12 &272 &   1 & 272&  23&  15 \\
10 &   1 &  10&   3&   3 &98 &   1 &  98&  11&   9 &186 &   1 & 186&  51&  12 &274 &   1 & 274&  37&  15 \\
12 &   1 &  12&   5&   3 &100 &   1 & 100&  13&   9 &188 &   1 & 188&  15&  13 &276 &   1 & 276&  37&  15 \\
{\bf 14} & {\bf 1} & {\bf 14} & {\bf 5} & {\bf 3} &102 &   1 & 102&  13&   9 &190 &   1 & 190&  79&  12 &278 &   1 & 278&  37&  15 \\
16 &   1 &  16&   5&   4 &104 &   1 & 104&  23&   9 &192 &   1 & 192&  35&  12 &280 &   1 & 280&  43&  14 \\
18 &   1 &  18&   5&   4 &106 &   1 & 106&  23&   9 &194 &   1 & 194&  35&  12 &282 &   1 & 282&  37&  15 \\
20 &   1 &  20&   5&   4 &108 &   1 & 108&  41&   9 &196 &   1 & 196&  17&  13 &284 &   1 & 284&  61&  15 \\
22 &   1 &  22&   5&   5 &110 &   1 & 110&  11&  10 &198 &   1 & 198&  17&  13 &286 &   1 & 286&  21&  16 \\
{\bf 24} & {\bf 2} & {\bf 12} & {\bf 3} & {\bf 4} &112 &   1 & 112&  25&   9 &200 &   1 & 200&  17&  13 &288 &   1 & 288&  77&  15 \\
26 &   1 &  26&   5&   5 &114 &   1 & 114&  25&   9 &202 &   1 & 202&  17&  13 &290 &   1 & 290&  39&  15 \\
28 &   1 &  28&   7&   5 &116 &   1 & 116&  13&  10 &204 &   1 & 204&  37&  12 &292 &   1 & 292&  39&  15 \\
30 &   1 &  30&   7&   5 &118 &   1 & 118&  15&  10 &206 &   1 & 206&  19&  13 &{\bf 294} & {\bf 7} & {\bf 42}& {\bf 3}& {\bf 14} \\
32 &   1 &  32&   7&   5 &120 &   1 & 120&  27&  10 &208 &   1 & 208&  57&  13 &296 &   1 & 296&  39&  15 \\
 34 &   1 &  34&   7&   5 &{\bf 122} & {\bf 1} & {\bf 122}& {\bf 27}& {\bf 9} &210 &   1 & 210&  77&  13 &298 &   1 & 298&  39&  16 \\
36 &   1 &  36&   7&   6 &124 &   1 & 124&  27&  10 &212 &   1 & 212&  33&  13 &300 &   1 & 300&  69&  15 \\
{\bf 38} & {\bf 1} & {\bf 38} & {\bf 15} & {\bf 5} &126 &   3 &  42&   5&  10 &214 &   1 & 214&  33&  13 &302 &   1 & 302&  65&  15 \\
40 &   1 &  40&   7&   6 &128 &   2 &  64&  13&  10 &{\bf 216} & {\bf 6} &  {\bf 36} & {\bf 3}& {\bf 12} &304 &   1 & 304&  41&  16 \\
42 &   1 &  42&   7&   6 &130 &   1 & 130&  29&  10 &218 &   1 & 218&  33&  13 &306 &   1 & 306&  41&  15 \\
44 &   1 &  44&   9&   6 &132 &   1 & 132&  29&  10 &220 &   1 & 220&  17&  14 &308 &   1 & 308&  41&  15 \\
46 &   1 &  46&   7&   7 &134 &   1 & 134&  13&  11 &222 &   1 & 222&  51&  13 &310 &   1 & 310&  41&  15 \\
48 &   1 &  48&  19&   6 &136 &   1 & 136&  13&  11 &224 &   4 &  56&   5&  13 &312 &   1 & 312&  67&  15 \\
50 &   1 &  50&   7&   7 &138 &   1 & 138&  13&  11 &226 &   1 & 226&  35&  13 &314 &   1 & 314&  59&  16 \\
52 &   1 &  52&   9&   7 &140 &   1 & 140&  31&  10 &228 &   1 & 228&  35&  13 &316 &   1 & 316&  59&  16 \\
{\bf 54} & {\bf 3} & {\bf 18} & {\bf 3} & {\bf 6} &142 &   1 & 142&  13&  11 &230 &   1 & 230&  35&  13 &318 &   1 & 318&  49&  16 \\
56 &   1 &  56&   9&   7 &144 &   1 & 144&  15&  11 &232 &   1 & 232&  89&  13 &320 &   1 & 320&  43&  16 \\
58 &   1 &  58&   9&   7 &146 &   1 & 146&  15&  11 &234 &   3 &  78&   7&  14 &322 &   1 & 322&  43&  15 \\
60 &   1 &  60&   9&   7 &148 &   1 & 148&  15&  11 &236 &   1 & 236&  43&  14 &324 &   1 & 324&  43&  15 \\
62 &   1 &  62&   9&   7 &{\bf 150} & {\bf 5} & {\bf 30}& {\bf 3}& {\bf 10} &238 &   1 & 238&  37&  14 &326 &   1 & 326&  43&  16 \\
64 &   1 &  64&  11&   7 &152 &   1 & 152&  27&  11 &240 &   1 & 240&  37&  13 &328 &   2 & 164&  23&  16 \\
66 &   1 &  66&  19&   7 &154 &   1 & 154&  13&  12 &242 &   1 & 242&  37&  13 &330 &   1 & 330&  71&  16 \\
68 &   1 &  68&  19&   7 &156 &   1 & 156&  43&  11 &244 &   1 & 244&  37&  14 &332 &   1 & 332&  51&  16 \\
70 &   1 &  70&   9&   8 &158 &   1 & 158&  29&  11 &246 &   1 & 246&  53&  14 &334 &   1 & 334&  21&  17 \\
72 &   1 &  72&   9&   8 &160 &   1 & 160&  29&  11 &248 &   1 & 248&  45&  14 &336 &   1 & 336&  45&  16 \\
{\bf 74} & {\bf 1} & {\bf 74} & {\bf 21} & {\bf 7} &162 &   1 & 162&  29&  11 &250 &   1 & 250&  17&  15 &{\bf 338} & {\bf 1} & {\bf 338}& {\bf 45}& {\bf 15} \\
76 &   1 &  76&  11&   8 &164 &   1 & 164&  45&  11 &252 &   1 & 252&  39&  14 &340 &   1 & 340&  45&  16 \\
78 &   1 &  78&  29&   8 &166 &   1 & 166&  17&  12 &{\bf 254} & {\bf 1} & {\bf 254}& {\bf 39}& {\bf 13} &342 &   3 & 114&  15&  16 \\
80 &   1 &  80&  23&   8 &168 &   1 & 168&  17&  12 &256 &   1 & 256&  39&  14 &344 &   1 & 344&  25&  17 \\
82 &   1 &  82&  23&   8 &170 &   1 & 170&  31&  11 &258 &   1 & 258&  99&  14 &346 &   1 & 346&  25&  17 \\
84 &   1 &  84&  11&   9 &172 &   1 & 172&  31&  11 &260 &   1 & 260&  19&  15 &348 &   1 & 348&  41&  17 \\
86 &   1 &  86&  11&   9 &174 &   1 & 174&  31&  12 &262 &   1 & 262&  19&  15 &350 &   1 & 350& 153&  16 \\
88 &   1 &  88&  25&   8 &176 &   1 & 176&  39&  12 &264 &   2 & 132&  19&  14 &352 &   1 & 352&  47&  16 \\
90 &   1 &  90&  11&   9 &178 &   1 & 178&  49&  12 &266 &   1 & 266&  41&  14 &354 &   1 & 354&  47&  16 \\
\hline
\end{tabular}
		\end{center}
		\end{table}

  \begin{table}[t]
  \caption{The minimum diameter $k$ for each number of vertices $N\le 526$ in the chordal ring mixed graphs $CRM(N,c)$ with (minimum) chord $c$. In boldface, there are the cases with a maximum number of vertices for a given diameter.}
		\label{tab:tabla-CRM}
  \small
		\begin{center}
\setlength{\tabcolsep}{3pt}
\begin{tabular}{|c|c|c||c|c|c||c|c|c||c|c|c||c|c|c||c|c|c|}
\hline
$N$ & $c$ & $k$ & $N$ & $c$ & $k$ & $N$ & $c$ & $k$ & $N$ & $c$ & $k$ & $N$ & $c$ & $k$ & $N$ & $c$ & $k$\\
\hline
{\bf 2} &   {\bf 1} &   {\bf 1} &	 90 &  33 &  13 &	178 &  17 &  19 &	266 &  23 &  23 &	354 &  75 &  27 &	442 & 139 &  29 \\
4 &   1 &   3 &	 92 &  11 &  15 &	180 &  17 &  19 &	268 &  61 &  23 &	356 &  23 &  27 &	444 &  27 &  31 \\
6 &   3 &   3 &	  94 &   11 &  13 &	182 &  19 &  19 &	270 &  97 &  23 &	358 &  23 &  27 &	446 &  27 &  29 \\
{\bf 8} &   {\bf 3} &   {\bf 3} &	 96 &   9 &  15 &	184 &  51 &  19 &	272 &  41 &  23 &	360 & 165 &  27 &	448 &  23 &  31 \\
{\bf 10} &   {\bf 3} &   {\bf 4} &	 {\bf 98} &  {\bf 13} &  {\bf 13} &	186 &  33 &  19 &	274 &  43 &  23 &	362 &  25 &  27 &	{\bf 450} &  {\bf 29} &  {\bf 29} \\
12 &   3 &   5 &	100 &  13 &  15 &	188 &  15 &  19 &	276 &  19 &  23 &	364 &  25 &  27 &	452 &  21 &  31 \\
14 &  3 &   5 &	102 &  39 &  15 &	190 &  41 &  19 &	278 &  65 &  23 &	366 &  27 &  27 &	454 &  81 &  31 \\
16 &   {\bf 3} &   {\bf 6} &	104 &   9 &  15 &	192 &  69 &  19 &	280 &  87 &  23 &	368 &  21 &  27 &	456 &  99 &  31 \\
 {\bf 18} &   {\bf 5} &   {\bf 5} &	106 &  11 &  15 &	194 &  85 &  19 &	282 & 127 &  23 &	370 & 115 &  27 &	458 &  53 &  31 \\
 20 &   3 &   7 &	108 &  33 &  15 &	196 &  17 &  19 &	284 &  21 &  23 &	372 &  87 &  27 &	460 &  55 &  31 \\
 22 &   5 &   7 &	110 &  13 &  15 &	198 &  17 &  21 &	286 &  21 &  25 &	374 &  49 &  27 &	462 &  25 &  31 \\
 24 &   5 &   7 &	112 &  13 &  15 &	{\bf 200} &  {\bf 19} &  {\bf 19} &	{\bf 288} &  {\bf 23} &  {\bf 23} &	376 &  51 &  27 &	464 &  25 &  31 \\
 26 &   7 &   7 &	114 &  15 &  15 &	{\bf 202} &  {\bf 91} &  {\bf 20} &	{\bf 290} & {\bf 133} &  {\bf 24} &	378 &  57 &  27 &	466 &  61 &  31 \\
 28 &   5 &   7 &	116 &  11 &  15 &	204 &  15 &  21 &	292 &  81 &  25 &	380 &  23 &  27 &	468 &  99 &  31 \\
 30 &   5 &   9 &	118 &  27 &  15 &	206 &  31 &  21 &	294 &  19 &  25 &	382 &  87 &  27 &	470 &  89 &  31 \\
 {\bf 32} &   {\bf 7} &   {\bf 7} &	120 &  33 &  15 &	208 &  37 &  21 &	296 &  19 &  25 &	384 & 117 &  27 &	472 &  23 &  31 \\
 {\bf 34} &  {\bf 13} &   {\bf 8} &	122 &  51 &  15 &	210 &  33 &  21 &	298 &  17 &  25 &	386 & 177 &  27 &	474 &  27 &  31 \\
 36 &  15 &   9 &	124 &  13 &  15 &	212 &  17 &  21 &	300 &  47 &  25 &	388 &  25 &  27 &	476 & 133 &  31 \\
 38 &   5 &   9 &	126 &  11 &  17 &	214 &  17 &  21 &	302 &  41 &  25 &	390 &  19 &  29 &	478 &  29 &  31 \\
 40 &   7 &   9 &	{\bf 128} &  {\bf 15} &  {\bf 15} &	216 &  99 &  21 &	304 &  21 &  25 &	{\bf 392} &  {\bf 27} &  {\bf 27} &	480 &  29 &  31 \\
 42 &   9 &   9 &	{\bf 130} &  {\bf 57} &  {\bf 16} &	218 &  15 &  21 &	306 &  21 &  25 &	{\bf 394} & {\bf 183} &  {\bf 28} &	482 &  31 &  31 \\
 {\bf 44} &  {\bf 13} &  {\bf 10} &	132 &  39 &  17 &	220 &  19 &  21 &	308 & 143 &  25 &	396 &  75 &  29 &	484 & 105 &  31 \\
 46 &   7 &   9 &	134 &  25 &  17 &	222 &  21 &  21 &	310 &  23 &  25 &	398 &  71 &  29 &	486 & 201 &  31 \\
 48 &   5 &  11 &	136 &  13 &  17 &	{\bf 224} &  {\bf 51} &  {\bf 22} &	312 &  23 &  25 &	400 &  91 &  29 &	488 &  25 &  31 \\
 {\bf 50} &   {\bf 9} &   {\bf 9} &	138 &  11 &  17 &	226 &  69 &  21 &	314 &  19 &  25 &	402 &  23 &  29 &	490 & 145 &  31 \\
 52 &   7 &  11 &	140 &  63 &  17 &	228 &  13 &  23 &	{\bf 316} &  {\bf 85} &  {\bf 26} &	404 &  23 &  29 &	492 &  57 &  31 \\
 54 &   7 &  11 &	142 &  15 &  17 &	230 &  17 &  21 &	318 &  69 &  25 &	406 &  47 &  29 &	494 &  59 &  31 \\
 56 &  21 &  11 &	144 &  15 &  17 &	232 &  15 &  23 &	320 &  37 &  27 &	408 &  93 &  29 &	496 &  65 &  31 \\
 58 &   9 &  11 &	146 &  17 &  17 &	234 &  69 &  21 &	322 & 133 &  25 &	410 &  21 &  29 &	498 &  87 &  31 \\
 60 &   7 &  11 &	{\bf 148} &  {\bf 41} &  {\bf 18} &	236 &  19 &  23 &	324 &  21 &  27 &	412 &  25 &  29 &	500 &  27 &  31 \\
 62 &  11 &  11 &	150 &  13 &  17 &	238 &  19 &  21 &	326 &  21 &  25 &	414 &  25 &  29 &	502 & 119 &  31 \\
 64 &  19 &  11 &	152 &  13 &  19 &	240 &  19 &  23 &	328 &  21 &  27 &	416 & 195 &  29 &	504 & 177 &  31 \\
 66 &  25 &  11 &	154 &  47 &  17 &	{\bf 242} &  {\bf 21} &  {\bf 21} &	330 & 117 &  25 &	418 &  27 &  29 &	506 & 235 &  31 \\
 68 &   9 &  11 &	156 &  15 &  19 &	244 &  21 &  23 &	332 &  17 &  27 &	420 &  27 &  29 &	508 &  29 &  31 \\
 70 &   9 &  13 &	158 &  15 &  17 &	246 &  17 &  23 &	334 &  23 &  25 &	422 &  29 &  29 &	510 &  29 &  33 \\
 {\bf 72} &  {\bf 11} &  {\bf 11} &	160 &  11 &  19 &	248 &  15 &  23 &	336 &  23 &  27 &	{\bf 424} &  {\bf 99} &  {\bf 30} &	512 &  31 &  31 \\
 {\bf 74} &  {\bf 31} &  {\bf 12} &	{\bf 162} &  {\bf 17} &  {\bf 17} &	250 &  67 &  23 &	{\bf 338} &  {\bf 25} &  {\bf 25} &	426 &  23 &  29 &	{\bf 512} &  {\bf 31} &  {\bf 31} \\
 76 &   9 &  13 &	164 &  13 &  19 &	252 &  47 &  23 &	340 &  25 &  27 &	428 &  19 &  31 &	{\bf 514} & {\bf 241} &  {\bf 32} \\
 78 &   9 &  13 &	166 &  31 &  19 &	254 &  45 &  23 &	342 & 123 &  27 &	430 &  91 &  29 &	516 & 135 &  33 \\
 80 &  35 &  13 &	168 &  45 &  19 &	256 &  19 &  23 &	344 &  45 &  27 &	432 &  21 &  31 &	518 &  23 &  33 \\
 82 &  11 &  13 &	170 &  39 &  19 &	258 &  19 &  23 &	346 &  21 &  27 &	434 & 189 &  29 &	520 & 165 &  33 \\
 84 &  11 &  13 &	172 &  15 &  19 &	260 &  55 &  23 &	348 &  21 &  27 &	436 &  25 &  31 &	522 &  93 &  33 \\
 86 &   9 &  13 &	174 &  15 &  19 &	262 &  21 &  23 &	350 &  55 &  27 &	438 &  25 &  29 &	524 &  71 &  33 \\
 {\bf 88} &  {\bf 19} &  {\bf 14} &	176 &  13 &  19 &	264 &  17 &  23 &	352 &  19 &  27 &	440 &  25 &  31 &	526 &  27 &  33 \\
\hline
\end{tabular}
		\end{center}
		\end{table}

\begin{table}[t]
\caption{The minimum diameter $k$ in the chordal multi-ring mixed graphs $CMRM(m,n,c)$ with minimum values of $m$ and $n$. In boldface, there are the cases with a maximum number of vertices for a given diameter.}
		\label{tab:tabla-cmrm}
  \small
		\begin{center}
\setlength{\tabcolsep}{3pt}
\begin{tabular}{|c|c|c|c|c||c|c|c|c|c||c|c|c|c|c||c|c|c|c|c|}
\hline
$N$ & $m$ & $n$ & $c$ & $k$ & $N$ & $m$ & $n$ & $c$ & $k$ & $N$ & $m$ & $n$ & $c$ & $k$ & $N$ & $m$ & $n$ & $c$ & $k$ \\
\hline
4 &   1 &   4&   3&   3 &92 &   2 &  46&   5&  13 &180 &   1 & 180&  17&  19 &268 &   1 & 268&  61&  23 \\
6 & 1 &   6 & 3 &   3 &94 &   1 &  94&  11&  13 &182 &   1 & 182&  19&  19 &270 &   1 & 270&  97&  23 \\
{\bf 8} & {\bf 1} & {\bf 8}& {\bf 3}& {\bf 3} &96 &   1 &  96&   9&  15 &184 &   1 & 184&  51&  19 &272 &   1 & 272&  41&  23 \\
{\bf 10} & {\bf 1} & {\bf 10}& {\bf 3}& {\bf 4} &{\bf 98} & {\bf 1} & {\bf 98}& {\bf 13}& {\bf 13} &186 &   1 & 186&  33&  19 &274 &   1 & 274&  43&  23 \\
12 &   1 &  12&   3&   5 &{\bf 100} & {\bf 2} & {\bf 50}& {\bf 7}& {\bf 14} & 188 &   1 & 188&  15&  19 &276 &   1 & 276&  19&  23 \\
14 &   1 &  14&   3&   5 &102 &   1 & 102&  39&  15 &190 &   1 & 190&  41&  19 &278 &   1 & 278&  65&  23 \\
16 &   1 &  16&   3&   6 &104 &   1 & 104&   9&  15 &192 &   1 & 192&  69&  19 &280 &   1 & 280&  87&  23 \\
{\bf 18} & {\bf 1} & {\bf 18}& {\bf 5}& {\bf 5} &106 &   1 & 106&  11&  15 &194 &   1 & 194&  85&  19 &282 &   1 & 282& 127&  23 \\
{\bf 20} & {\bf 2} & {\bf 10}& {\bf 3}& {\bf 6} &108 &   1 & 108&  33&  15 &196 &   1 & 196&  17&  19 &284 &   1 & 284&  21&  23 \\
22 &   1 &  22&   5&   7 &110 &   1 & 110&  13&  15 &198 &   9 &  22&  11&  19 &286 &  11 &  26&  13&  23 \\
24 &   1 &  24&   5&   7 &112 &   1 & 112&  13&  15 &{\bf 200} & {\bf 1} & {\bf 200}& {\bf 19}& {\bf 19} & {\bf 288} & {\bf 1} & {\bf 288}& {\bf 23}& {\bf 23} \\
 26 &   1 &  26&   7&   7 &114 &   1 & 114&  15&  15 &{\bf 202} & {\bf 1} & {\bf 202}&  {\bf 91}& {\bf 20} &{\bf 290} & {\bf 1} & {\bf 290}& {\bf 133}&  {\bf 24} \\
28 &   1 &  28&   5&   7 &116 &   1 & 116&  11&  15 &204 &   1 & 204&  15&  21 &{\bf 292} & {\bf 1} & 292& 81& 25 \\
30 &   3 &  10&   5&   7 &118 &   1 & 118&  27&  15 &206 &   1 & 206&  31&  21 &294 &   1 & 294&  19&  25 \\
{\bf 32} & {\bf 1} & {\bf 32} & {\bf 7} & {\bf 7} &120 &   1 & 120&  33&  15 &208 &   1 & 208&  37&  21 &296 &   1 & 296&  19&  25 \\
{\bf 34} & {\bf 1} & {\bf 34}& {\bf 13}& {\bf 8} &122 &   1 & 122&  51&  15 &210 &   1 & 210&  33&  21 &298 &   1 & 298&  17&  25 \\
36 &   1 &  36&  15&   9 &124 &   1 & 124&  13&  15 &212 &   1 & 212&  17&  21 &300 &   1 & 300&  47&  25 \\
38 &   1 &  38&   5&   9 &126 &   7 &  18&   9&  15 &214 &   1 & 214&  17&  21 &302 &   1 & 302&  41&  25 \\
40 &   1 &  40&   7&   9 &{\bf 128} & {\bf 1} & {\bf 128}& {\bf 15}& {\bf 15} &216 &   1 & 216&  99&  21 &304 &   1 & 304&  21&  25 \\
42 &   1 &  42&   9&   9 &{\bf 130} & {\bf 1} & {\bf 130}& {\bf 57}& {\bf 16} &218 &   1 & 218&  15&  21 &306 &   1 & 306&  21&  25 \\
44 &   2 &  22&   3&   9 &132 &   1 & 132&  39&  17 &220 &   1 & 220&  19&  21 &308 &   1 & 308& 143&  25 \\
46 &   1 &  46&   7&   9 &134 &   1 & 134&  25&  17 &222 &   1 & 222&  21&  21 &310 &   1 & 310&  23&  25 \\
48 &   1 &  48&   5&  11 &136 &   1 & 136&  13&  17 &224 &   1 & 224&  51&  22 &312 &   1 & 312&  23&  25 \\
{\bf 50} & {\bf 1} & {\bf 50}& {\bf 9}& {\bf 9} &138 &   1 & 138&  11&  17 &226 &   1 & 226&  69&  21 &314 &   1 & 314&  19&  25 \\
{\bf 52} & {\bf 2} & {\bf 26}& {\bf 5}& {\bf 10} &140 &   1 & 140&  63&  17 &228 &   2 & 114&  33&  21 &316 &   2 & 158&  33&  25 \\
54 &   1 &  54&   7&  11 &142 &   1 & 142&  15&  17 &230 &   1 & 230&  17&  21 &318 &   1 & 318&  69&  25 \\
56 &   1 &  56&  21&  11 &144 &   1 & 144&  15&  17 &232 &   4 &  58&  25&  21 &320 &   1 & 320&  37&  27 \\
58 &   1 &  58&   9&  11 &146 &   1 & 146&  17&  17 &234 &   1 & 234&  69&  21 &322 &   1 & 322& 133&  25 \\
60 &   1 &  60&   7&  11 &148 &   2 &  74&  23&  17 &236 &   2 & 118&   9&  21 &324 &   2 & 162&  57&  25 \\
62 &   1 &  62&  11&  11 &150 &   1 & 150&  13&  17 &238 &   1 & 238&  19&  21 &326 &   1 & 326&  21&  25 \\
64 &   1 &  64&  19&  11 &152 &   4 &  38&  13&  17 &240 &   1 & 240&  19&  23 &328 &   4 &  82&  33&  25 \\
66 &   1 &  66&  25&  11 &154 &   1 & 154&  47&  17 &{\bf 242} & {\bf 1} & {\bf 242}& {\bf 21}& {\bf 21} &330 &   1 & 330& 117&  25 \\
68 &   1 &  68&   9&  11 &156 &   2 &  78&   7&  17 &{\bf 244} & {\bf 2} & {\bf 122}& {\bf 11}& {\bf 22} &332 &   2 & 166&  11&  25 \\
70 &   5 &  14&   7&  11 &158 &   1 & 158&  15&  17 &246 &   1 & 246&  17&  23 &334 &   1 & 334&  23&  25 \\
{\bf 72} & {\bf 1} & {\bf 72} & {\bf 11} & {\bf  11} &160 &   1 & 160&  11&  19 &248 &   1 & 248&  15&  23 &336 &   1 & 336&  23&  27 \\
 {\bf 74} & {\bf 1} & {\bf 74}& {\bf 31}& {\bf 12} &{\bf 162} & {\bf 1} & {\bf 162}& {\bf  17}& {\bf 17} &250 &   1 & 250&  67&  23 &{\bf 338} & {\bf 1} & {\bf 338}& {\bf 25}& {\bf 25} \\
76 &   1 &  76&   9&  13 &{\bf 164} & {\bf 2} & {\bf 82}& {\bf 9}& {\bf 18} &252 &   1 & 252&  47&  23 &{\bf 340} & {\bf 2} & {\bf 170}& {\bf 13}& {\bf  26} \\
78 &   1 &  78&   9&  13 &166 &   1 & 166&  31&  19 &254 &   1 & 254&  45&  23 &342 &   1 & 342& 123&  27 \\
80 &   1 &  80&  35&  13 &168 &   1 & 168&  45&  19 &256 &   1 & 256&  19&  23 &344 &   1 & 344&  45&  27 \\
82 &   1 &  82&  11&  13 &170 &   1 & 170&  39&  19 &258 &   1 & 258&  19&  23 &346 &   1 & 346&  21&  27 \\
84 &   1 &  84&  11&  13 &172 &   1 & 172&  15&  19 &260 &   1 & 260&  55&  23 &348 &   1 & 348&  21&  27 \\
86 &   1 &  86&   9&  13 &174 &   1 & 174&  15&  19 &262 &   1 & 262&  21&  23 &350 &   1 & 350&  55&  27 \\
88 &   4 &  22&   9&  13 &176 &   1 & 176&  13&  19 &264 &   1 & 264&  17&  23 &352 &   1 & 352&  19&  27 \\
90 &   1 &  90&  33&  13 &178 &   1 & 178&  17&  19 &266 &   1 & 266&  23&  23 &354 &   1 & 354&  75&  27 \\
\hline
\end{tabular}
		\end{center}
		\end{table}

\end{document}